\documentclass[twoside]{article}
\usepackage{eurosym}
\usepackage[utf8]{inputenc}
\usepackage{graphicx}
\usepackage{latexsym}
\usepackage{amsmath}
\usepackage{amsfonts}
\usepackage{amssymb}
\usepackage{hyperref}
\usepackage[usenames]{color}

\newcommand{\disp}{\displaystyle}
\newcommand{\1}{1\!{\rm I}}
\newcommand{\cA}{{\mathcal A}}
\newcommand{\rA}{{\rm A}}

\newcommand{\rB}{{\rm B}}
\newcommand{\bB}{{\bf B}}
\newcommand{\cB}{{\mathcal B}}
\newcommand{\BB}{\mathbb{B}}
\newcommand{\cC}{{\mathcal C}}
\newcommand{\rC}{{\rm C}}
\newcommand{\rc}{{\rm c}}
\newcommand{\rD}{{\rm D}}
\newcommand{\cD}{{\mathcal D}}
\newcommand{\rE}{{\rm E}}
\newcommand{\re}{{\rm e}}
\newcommand{\EE}{\mathbb{E}}

\newcommand{\rF}{{\rm F}}
\newcommand{\rf}{{\rm f}}

\newcommand{\cF}{{\mathcal F}}

\newcommand{\cG}{{\mathcal G}}
\newcommand{\HH}{\mathbb{H}}
\newcommand{\rH}{{\rm H}}
\newcommand{\bH}{{\bf H}}

\newcommand{\rI}{{\rm I}}

\newcommand{\rK}{{\rm K}}

\newcommand{\cL}{{\mathcal L}}
\newcommand{\rL}{{\rm L}}

\newcommand{\rM}{{\rm M}}

\newcommand{\cN}{{\mathcal{N}}}

\newcommand{\PP}{{\mathbb{P}}}
\newcommand{\cP}{{\mathcal P}}
\newcommand{\rP}{{\rm P}}
\newcommand{\cQ}{{\mathcal Q}}
\newcommand{\cQh}{\widehat{\mathcal Q}}
\newcommand{\rQ}{{\rm Q}}
\newcommand{\RR}{{\mathbb{R}}}

\newcommand{\rR}{{\rm R}}
\newcommand{\rS}{{\rm S}}

\newcommand{\rT}{{\rm T}}

\newcommand{\rU}{{\rm U}}
\newcommand{\ru}{{\rm u}}

\newcommand{\rV}{{\rm V}}

\newcommand{\bv}{{\bf v}}
\newcommand{\rv}{{\rm v}}

\newcommand{\rW}{{\rm W}}

\newcommand{\rw}{{\rm w}}
\newcommand{\rX}{{\rm X}}

\newcommand{\rY}{{\rm Y}}

\newcommand{\fin}{\hfill\mbox{$\quad{}_{\Box}$}}
\newcommand{\fineq}{\vspace{-.75cm$\fin$}\par\bigskip}
\newcommand{\fineqnum}{\vspace{-.4cm$\fin$}\par\bigskip}
\newcommand{\pe}[2]{\langle #1,#2\rangle}
\newcommand{\n}[1]{\disp \|#1\|}
\newcommand{\nl}[1]{\disp \big \|#1\big \|}
\newcommand{\nll}[1]{\disp \bigg \|#1\bigg \|}
\newcommand{\nn}[1]{\disp \|\!|#1\|\!|}

\newtheorem{theo}{\bf \sffamily Theorem}
\newtheorem{prop}{\bf \sffamily Proposition}

\newtheorem{rem}{\bf \sffamily Remark}
\newtheorem{lemma}{\bf \sffamily Lemma}

\begin{document}

\title{\Large \bfseries\sffamily Stochastic diffusive energy balance climate model with a multiplicative noise modeling the Solar variability}
\author{\bfseries\sffamily G. D\'{\i}az \& J.I. D\'{\i}az\thanks {The research of J.I. Díaz was partially supported by the project PID2020-112517 GB-I00 of the Spain State Research Agency (AEI) and PID2023-146754NB-I00 funded by MCIU/AEI/10.13039/501100011033 and FEDER, EU. MCIU/AEI/10.13039/-501100011033/FEDER,
EU. \hfil\break \indent {\sc
			Keywords}: Stochastic Diffusive Energy Balance Models, parabolic Legendre diffusion, multiplicative noise, cylindrical Wiener process, continuous non-Lipschitz hybrid co-albedo, comparison of solutions, successive approximations.\hfil\break \indent {\sc AMS Subject Classifications: 86A08, 60H15, 35K55, 35Q79. 60H30.}}}
\date{\boxed{\hbox{On September 24, 2025}}}
\maketitle

\begin{center}
Dedicated to Roger Temam, always admired, on the occasion of his 85th
birthday
\end{center}
\begin{abstract}
We prove existence, uniqueness, and comparison of
solutions for a nonlinear stochastic parabolic partial
differential equation that includes the Solar variability in terms of a
multiplicative Wiener cylindrical noise in the term of the absorbed
radiative energy in a simplified diffusive one-dimensional Energy Balance
Model. We introduce a hybrid co-albedo nonlinear term, which has the
advantages of both the Sellers model, as it is a continuous function, and the
Budyko model, as it has an infinite derivative at $u=-10^{\circ}C$ (the temperature
at which ice is white), allowing the location of the polar ice caps to be
easily detected. We show that, despite the lack of differentiability of this function, the method of successive approximations can be satisfactorily applied.
\end{abstract}

\section{Introduction}
The radiative energy balance climate models (EBMs) are a class of tools for
representing the evolution of the global climate (spatial distribution of
temperature at the Earth's surface) over large time scales. Despite their
simplicity, the EBMs give a useful, representation of Earth's climate by
capturing the fundamental mechanisms governing its behaviour. They were
proposed in 1969, independently (during the "Cold War") by an American,
Sellers \cite{sellers}, and a Russian, Budyko \cite{Budiko}, maintaining an
enormous resemblance that confirms their great robustness. The EBMs assume  
that the averaged atmospheric temperature evolves (on a large scale of time)
according to the radiation balance of the budget, i.e. the difference in the
radiations absorbed and emitted by the planet. The Solar radiation (is the
primary input) though must be corrected by the Earth's co-albedo (which depends
of the temperature as a feedback) and the balance is taken with the infrared
radiation emitted by the Earth and also in the presence of a surface
diffusion. When describing the constitutive laws of the absorbed and emitted
radiation it is needed to take into account important elements such as the
Solar constant, the Solar insolation (depending of the spatial
distribution), the atmospheric composition (which appears as coefficient
when applying the Stefan-Boltzmann for the outgoing longwave radiation) etc.
Since the pioneering work of 1969, mentioned above, EBMs have been the
subject of numerous research studies and monographs (\cite{NorthLibro}, \cite%
{Ghil-Childres}, \cite{Hetzer1990}, \cite{D-Escorial}, \cite{Henerson-sellers}, 
\cite{Ghil-Lucarini}, \cite{Arcoya-D-T}, \cite{DHetzTello}, \cite%
{Bermejo-Diaz}, \cite{D-Hetzer}, \cite{diazNATO}, \cite{Kaper}, \cite%
{Bensid-D}, \cite{Canarsa-Lucarini}, etc.). Some version of the EBMs can be
also obtained by averaging in the {\em primitive equations} (see, e.g., \cite%
{Lions-Temam-Wang92}, \cite{LTM92b}, \cite{LTW95}), as it was presented in 
\cite{Kiel}. A different averaging approach was developed by Hasselmann 
\cite{Hasselman} (see also \cite{Arnold} and \cite{del2024non}).
EBMs are also very useful in the study of past climates (\cite{CrowleyNorth}%
).

Deterministic EBMs do not include unpredictable external forcings, as, for
instance, volcanic emissions (there are currently more than 500 active
volcanoes), etc. From the mathematical point of view, this has been treated by
means of an additive white noise (\cite{North-Calahan}, \cite{GD-JID2002}, 
\cite{DLangaValero}, \cite{Imkeller}, \cite{GD-JID2022}, \cite{Lucarini}, 
\cite{delSartoBroknerFlandoli}, \cite{del2024non}). The
mathematical models can be also coupled with some simple modelling of the
deep ocean temperature (\cite{Watts-Morantine}, \cite{DTello}, \cite%
{D-TelloOcean}, \cite{Diaz-Hidalgo-Tello}, \cite{delasarto2025}) but, for
simplicity in the formulation, we will not follow this coupling in this
paper.

EBMs with stochastic noise allow us to justify, through scientific
arguments and the available data, the possible increase in extreme events
due to Climate Change, that is not possible under purely deterministic
approaches (see, e.g. \cite{del2024non}).

The main goal of this paper is to study a mathematical model taking into
account the influence on the climate of the abrupt changes in the Solar
radiation (the Solar storms). The assumption that Solar emission is constant
must be replaced by a more realistic study that takes Solar variability into
account (\cite{Scafetta-survey}, \cite{Vaquero-VazquezLibro}). In fact,
there is a whole family of space satellites whose primary mission is to
analyze and measure this Solar variability. The Earth Radiation Budget
Satellite (ERBS) was a NASA scientific research satellite. The satellite was
one of three satellites in NASA's research program, named Earth Radiation
Budget Experiment (ERBE), to investigate the Earth's radiation budget.
NASA's CERES instruments have continued the ERB data record after 1997. We
recall that the Solar energy that falls annually on Earth's
surface is about ten thousand times the energy demand of the world's
population (7.7 billion people). Phenomena such as Solar flares or coronal
mass ejections (extreme Solar events) can cause brief increases in
radiation. Occasionally, flares heat the Sun's surface, reaching
temperatures of about 45 million degrees Celsius, much higher than those in
the core. Long-term variations are minimal (the magnitude of the short-term
fluctuations are small, typically less than 0.1\% of the average value) but
relevant in the climate and paleoclimate studies (\cite{CrowleyNorth}).

The main factor that generates variations is the Solar activity. The
so-called\ S. Schabe (1789-1875) cycles vary over 11 years (\cite{MS}). This is too short a period to have any impact on the climate,
though it is very relevant in other aspects. This was previously measured by
counting sunspots, and is now measured using satellite radiance
measurements. Other types of cycles (the so-called W. Gleissberg (1903-1986)
cycle)) has an oscillation amplitude similar to the Schwabe cycles but its
duration is approximately 87 years (70-100 years) and has a greater impact
on climate due to its duration (\cite{peristyk}). It is related to the well-identified past periods (Maunder Minimum (1645-1715) and the Dalton Minimum
(1800-1830)) of an extraordinarily low Solar activity. There is also the
well-known Milankowitz cycles, based on Celestial Mechanics, which are
justified in another way (on scales of thousands of years) and can be
considered as periodic versions of the Solar constant (see a mathematical
study\ on a pure time periodical Solar datum $Q(t)$ in \cite%
{Badii}).
\par
\medskip
The mathematical model that we will consider in this paper includes the
Solar variability in terms of a multiplicative Wiener noise in the term of
the absorbed radiative energy in a simplified diffusive one-dimensional
Energy Balance Model:

\begin{equation*}
(\mathrm{E}_{\beta ,\varepsilon })\left\{ 
\begin{array}{l}
du_{t}-\dfrac{\partial }{\partial x}\left( \big (1-x^{2})\dfrac{\partial
	u_{t}}{\partial x}\right) +g(u_{t})=\mathrm{Q}\mathrm{S}(x)\beta
(u_{t})\left( 1+\varepsilon d\mathrm{W}_{t}\right) , \\[0.2cm]
u(x,0)=u_{0}(x).%
\end{array}%
\right. 
\end{equation*}%
where $x\in \mathrm{I}\doteq (-1,1)$, $x=\cos \phi ,$ with $\phi $ the
spherical latitude,  $t>0$ and $\varepsilon \geq 0$. Notice that we are
following the usual dynamical system notation (see, e.g., \cite%
{TemamLibro}), $u=u(x,t;\omega )=u_{t}(x;\omega )$ where $x\in \mathrm{I}%
,\doteq (-1,1),$ $x=\cos \phi ,$ with $\phi $ the spherical latitude, $%
~t\geq 0$ and $\omega $ is in the probability space $\{\Omega ,{\mathcal{F}},{%
\mathbb{P}}\}$. The {\em cylindrical Wiener processes} $\rW_{t}(x;\omega )$ 
is not time differentiable and thus the notation
used for deterministic models $\frac{\partial u}{\partial t}$\ is not
well justified. In what follows, we will use the notation $\rho
(x)=1-x^{2},~x\in \mathrm{I}$ for the degenerate diffusion operator
coefficient. Note that in this formulation, it is not necessary to specify any boundary
condition on $\partial \mathrm{I}$, since the physical problem is posed on
the sphere as a Riemannian manifold without boundary. This explains the
degeneracy of the boundary operator on $\partial \mathrm{I}$.
\par
\medskip
Although we can assume greater generality, here we mainly assume that $\mathrm{Q}>0$, and that the Earth emitted radiation is given by the term $g(u)$ such that
\par
(\textbf{H}$_{g}$) $g:{\mathbb{R}}\rightarrow {\mathbb{R}}$ is a continuous
increasing function.
\par
\noindent The co-albedo term is assumed to be such that
\par
(\textbf{H}$_{\beta }$) $\beta $ is a \emph{bounded} maximal monotone graph
in ${\mathbb{R}}^{2},~$i.e.$~\beta (s)\in \lbrack m,\mathrm{M}]$, $\forall
s\in {\mathbb{R}}$.
\par
\medskip
We recall that one of the important differences in the 1969
modelling of these terms was concerning their regularity: In  \cite%
{sellers} it was assumed that $\beta (s)$ is a Lipschitz continuous
non-decreasing function with a huge slope near $u=-10^{\circ}$ (the Celsius
temperature at which the ice goes from transparent to white), in contrast
with \cite{Budiko}, who assumed that the co-albedo is a discontinuous
function at $u=-10,$ to better parameterize the regions occupied by
the polar caps. Here we are identifying the possible discontinuous
non-decreasing function $\beta (s)$ with the associated maximal monotone 
graph (\cite{brezis1973ope}) by including the entire jump interval at the discontinuity points.
Although the result leads to a possible multivalued expression, we will
simplify the writing by keeping the symbol $=$ in the equation.
\par
\medskip
On the \textit{insolation function} $\mathrm{S}(x)$ we assume
\par
(\textbf{H}$_{s}$) $\mathrm{S}:\mathrm{I}\rightarrow {\mathbb{R}},~\mathrm{S}%
\in \mathrm{L}^{\infty }(\mathrm{I})$, $\mathrm{S}_{1}\geq \mathrm{S}(x)\geq 
\mathrm{S}_{0}>0$ $\ \mbox{a.e.}~x\in \mathrm{I}$.
\par
\begin{figure}[htp]
\begin{center}
\vspace*{.5cm}
\includegraphics[width=10cm]{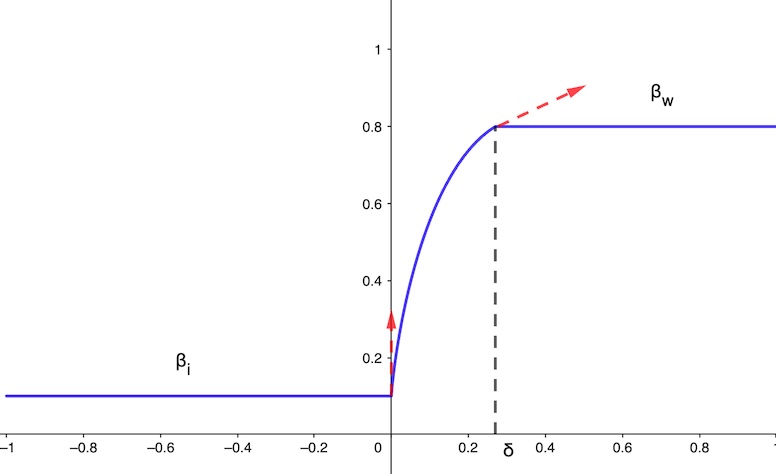}\\ 
\caption{Hybrid co-albedo profile}
\label{fig:coalbedoprofile}
\end{center}
\end{figure}

Despite the presence of multiplicative noise $d\mathrm{W}_{t}$, there are
many abstract results which can be applied when $\beta $ is, as in the
Sellers' option, a globally Lipschitz continuous function of the unknown
(see, e.g., \cite{da2014stochastic}). On the other hand, the study of the
multiplicative noise in the presence of a discontinuous co-albedo function (as
in the Budyko case) looks very complex to be considered in a first
approach. In that paper, we will follow an intermediate option by considering
a class of co-albedo functions $\beta $ that, being continuous, present a
singularity in their derivative at the critical value $u=-10$, which,
following Budyko's motivation, allows us to easily recognize the regions
occupied by the polar caps. 
More specifically, let us define the co-albedo function given by 
\begin{equation}
\beta _{-10}(u)=\left\{ 
\begin{array}{ll}
\beta _{i}, & \quad \hbox{if $u<-10$,} \\[0.1cm]
\big (\beta _{w}-\beta _{i}\big )\theta _{\delta +10}(u+10)+\beta _{i}, & 
\quad \hbox{if $-10\le u\le -10+\delta $}, \\[0.1cm]
\beta _{w}, & \quad \hbox{if $u>-10+\delta$},%
\end{array}%
\right.   \label{eq:coalbedofunction}
\end{equation}%
with $0<\delta <1,$ where the function $\theta _{\delta }(u)$ is given by
\begin{equation}
\theta _{\delta }(u)=(\beta_{w}-\beta _{i})\dfrac{u \ln u}{\delta \ln \delta},\quad u\ge 0.
\label{eq:thetaintro}
\end{equation}
We note that $\beta \in {\mathcal{C}}({\mathbb{R}})\cap {%
\mathcal{C}}^{1}\big ({\mathbb{R}}\setminus \{-10,-10+\delta \}\big )$, with 
\begin{equation*}
\left\{ 
\begin{array}{l}
\beta _{-10}^{\prime }(-10^{-})=0\quad \hbox{and \quad }\beta _{-10}^{\prime
}(-10^{+})=+\infty , \\ 
\beta _{-10}^{\prime }\big ((-10+\delta )^{-}\big )>0\quad \hbox{and \quad }%
\beta _{-10}^{\prime }\big ((-10+\delta \big )^{+})=0.%
\end{array}%
\right. 
\end{equation*}%
Thus, this choice of the co-albedo function (that we will call the {\em hybrid co-albedo function} in what follows) presents a sudden change at the critical temperature for which ice becomes white while offering a seamless transition from ice to water (see Figure \ref{fig:coalbedoprofile} below where the profile of the co-albedo function was transferred to the origin).
\par
\medskip
We recall that the uniqueness of solutions when the right-hand side term is not Lipschitz continuous is a very delicate question. For deterministic parabolic PDEs, an important contribution was offered in the paper \cite{fujita1968uniqueness}. They prove that the usual uniqueness criterion for non-monotone elliptic equations is not enough for the parabolic equation, and they prove the uniqueness of solutions by asking an Osgood-type condition on the non-linear term. We will extend this type of results, in the stochastic framework, in two different directions which, in the best of our knowledge, were not presented in \cite{fujita1968uniqueness}, nor
on its generalizations): our proof will be constructive (since we will
prove that the successive approximations method can be applied even in the
absence of differentiability on $\beta $). In addition, we will get some
comparison results (in terms of two different initial data).
\par
\medskip
Faced with a wide range of possible choices, in this article, we will model
the  \emph{cylindrical Wiener noise }$\mathrm{W}_{t}(x;\omega )$ produced by
the erratic Solar storms using a series expansion of the
eigenfunctions generated by the Legendre diffusion mentioned above. We
recall that  the sequence of eigenvalues are $\mu_{n}=n(n+1)$ and that after a normalization the eigenfunctions are given by 
$$
\re_{n}(x)=\sqrt{\dfrac{2}{2n+1}}\rP_{n}(x),\quad -1\le x\le 1,\\ [.3cm]
$$
with (the Rodrigues formula)
$$
\rP_{n}(x)=\dfrac{1}{2^{n}n!}\dfrac{d^{n}}{dx ^{n}}\big (x^{2}-1\big )^{n}=\dfrac{1}{2^{n}}\sum_{k=0}^{n}\binom{n}{k}^{2}(x+1)^{n-k}(x-1).
$$
Since $\mu _{0}=0$, to avoid strong difficulties, we artificially
introduce a given (but arbitrary) $\mu >0$ and then we will work with the
perturbed diffusion $-\dfrac{\partial }{\partial x}\left( \big (1-x^{2})%
\dfrac{\partial u}{\partial x}\right) +\mu u$ (the artificial term $\mu u
$ will also be added to the corresponding right-hand side). So that the same
family of eigenfunctions corresponds now to the new sequence of eigenvalues  $ \widehat{\mu}_{n}=\mu_{n}+\mu>0$. Thus, the noise we consider will be given by 
\begin{equation}
\rW_{t}=\sum_{n\ge 0}\dfrac{1}{\sqrt{\widehat{\mu}}_{n}}\bB_{t}^{n}\re_{n},\quad t\ge 0,
\label{eq:WienerQIntro}
\end{equation}
where the processes $\{\bB^{n}_{t}\}_{t\ge 0}$ is a family of Brownian motions mutually independent (see \cite{chow2007stochastic,da2014stochastic,liu2015stochastic}, and Section 3, for a definition). We note that finally the parameter $\mu$ only appear in the noise (see~\eqref{eq:WienerQIntro}).

It is important to highlight that the main difficulty in the study of the multiplicative noise comes from the presence of the non-linear term $\beta(u)$ since if this term were a linear function in $\beta(u)$ the problem would be reduced to applying a clever change of variable that leads to a deterministic problem dependent on a parameter (see Remark \ref{Rmchangevariabl}).

The main contribution in this paper is 
\begin{theo} Let us assume $\ru_{0}\in \rL^{\infty}(\rI)$, let $\mu >0$ arbitrarily given. Assume {\rm (\textbf{H}$_{g}$), (\textbf{H}$_{s}$)} and $\beta $ given by (\ref%
{eq:coalbedofunction}). Then there exists a unique mild solution $\ru^{\ru_{0}}$, in the sense of \eqref{eq:mildsolutionmuFB}, of the climate diffusive energy balance model $(\rE_{\beta ,\varepsilon})$.
Furthermore, if $\ru^{\ru_{0}}$ and $\ru^{\widehat{\ru}_{0}}$ denote the solutions of the  problems relative to these data, one has the continuous dependence inequality 
\begin{equation}
\disp \nl{\ru^{\ru_{0}}-\ru^{\widehat{\ru}_{0}}}_{\BB_{t}}^{2}\le
4\rM^{2}\n{\ru_{0}-\widehat{\ru}_{0}}_{\rH}^{2}+\widehat{\rC}_{\rT} \int^{\rT}_{0}\theta_{\rF,\rB}\big (\nl{\ru^{\ru_{0}}-\ru^{\widehat{\ru}_{0}}}_{\BB_{s}}^{2}\big )ds 
\label{eq:dependenceu0fIntro}
\end{equation}
where the constants $\rM,\widehat{\rC}_{\rT}$ and the function $\theta_{\rF,\rB}$ are given  in Proposition \ref{prop:Gcontinuous}  and the Banach  space $\BB_{t}$ in \eqref{eq:BTspace}. In addition, introducing a suitable function \eqref{eq:Psifunction}, we have the growth estimate
\begin{equation}
\Psi_{4\rM^{2}\n{\ru_{0}-\widehat{\ru}_{0}}_{\rH}^{2}} \left ( \n{\ru^{\ru_{0}}-\ru^{\widehat{\ru}_{0}}}_{\BB_{t}}^{2}\right )\le \widehat{\rC}_{\rT}t, \quad t\in [0,\rT],
\label{eq:integraldependenceu0fintro}
\end{equation}
provided $\n{\ru_{0}-\widehat{\ru}_{0}}_{\rH}>0$. Finally, we have the quantitative comparison estimate
\begin{equation}
\disp \nl{\big(\ru^{\ru_{0}}-\ru^{\widehat{\ru}_{0}}\big )_{+}}_{\BB_{t}}^{2}\le
4\rM^{2}\n{\big (\ru_{0}-\widehat{\ru}_{0}\big )_{+}}_{\rH}^{2}+\widehat{\rC}_{\rT} \int^{\rT}_{0}\theta_{\rF,\rB}\big (\nl{\big(\ru^{\ru_{0}}-\ru^{\widehat{\ru}_{0}}\big )_{+}}_{\BB_{s}}^{2}\big )ds 
\label{eq:comparisonu0fintro}
\end{equation}
for $t\in [0,\rT]$ holds, where $r_{+}=\max \{r,0\}$. 
In consequence, we have the comparison of solutions: 
$$
\ru_{0}\le \widehat{\ru}_{0}\quad \hbox{in $\rL^{2}(\rI)$}
\quad \Rightarrow \quad \ru^{\ru_{0}}_{t}(x;\omega)\le \ru^{\widehat{\ru}_{0}}_{t}(x;\omega)
\quad \hbox{for all $(x,t)\in\rI\times [0,\rT],~a.e.~\omega\in\Omega$}.
$$
\label{theo:mainTheoremintro}
\end{theo}
\vspace*{-.3cm}
We point out that in a separate study (\cite{diaz25deterministic})  we will offer some extension and related results to those presented in this article, but in some different contexts: we will prove a strictly deterministic version of the aforementioned improvements of the important article \cite{fujita1968uniqueness}, and we will also explain how the stochastic problem can be addressed by applying the results in an abstract framework that offers possible applications to problems not necessarily originating from climate models (\cite{diaz25abstract}).

The organization of this paper is the following. Some useful properties of the above-presented hybrid co-albedo function will be given in Section 2. In Section 3 we collect some auxiliary and technical results which will be useful in this stochastic framework. Finally, the proof of the main Theorem  will be given in Section 4.

\section{On the hybrid co-albedo function}
\label{sec:coalbedofunction}
In this Section, we collect some results on some properties of the concave and increasing continuous function $\theta _{\delta }(u)$ used in the definition of the hybrid co-albedo function. For simplicity in the notation, we drop the dependence with respect to the constant $\delta$. So, let
\begin{equation}
\theta (u)=(\beta_{w}-\beta _{i})\dfrac{u \ln u}{\delta \ln \delta},\quad u\ge 0.
\label{eq:theta}
\end{equation}
For Section 4, it will be very useful the study of the nonlinear integral equation  
\begin{equation}
\rv(t)=\rv_{0}+\alpha \int^{t}_{0}\theta \big (\rv(s)\big )ds, \quad 0\le s\le \rT.
\label{eq:mainequation}
\end{equation}
It will have a key role in proving the existence part of the main theorem.  
\begin{theo} 
Let $\rv_{0},\alpha >0$ and let  $\rv(t)$ be any nonnegative integrable function satisfying 
\begin{equation}
\rv(t)\le \rv_{0}+\alpha \int^{t}_{0}\theta \big (\rv(s)\big )ds, \quad 0\le s\le \rT\le +\infty.
\label{eq:maininequation}
\end{equation}
Then we have the implicit estimate, 
\begin{equation}
\int_{\rv_{0}}^{\rv(t)}\dfrac{ds}{\theta (s)}\le \alpha t,\quad 0\le t\le \rT.
\label{eq:Leibnitzproperty}
\end{equation} 
In fact, the unique nonnegative function  $\rv(t)$ such that
\begin{equation}
\rv(t)\le \alpha \int^{t}_{0}\theta \big (\rv(s)\big )ds, \quad 0\le s\le \rT
\label{eq:nullinequation}
\end{equation}
is the null function. Moreover, if $\rv_{0}\ge 0$ and $\alpha >0$, under the Osgood condition
\begin{equation}
\int_{0^{+}}\dfrac{du}{\theta(u)}=\infty,
\label{eq:Osgoodcondition}
\end{equation}
the nonlinear integral equation 
\begin{equation}
\rv(t)= \rv_{0}+\alpha \int_{0}^{t}\theta \big (\rv(s)\big )ds,\quad 0\le t\le \rT,
\label{eq:a2}
\end{equation}
admits a unique global nonnegative solution, on $[0,\rT]$.
\label{theo:coalbedohybridfunction}
\end{theo}

{\sc Proof} As it is well known, if $\rv_{0}>0$, the integral equation \eqref{eq:mainequation} is equivalent to the Cauchy problem
\begin{equation}
\left \{
\begin{array}{l}
\rv'(t)=\alpha \theta \big (\rv(t)\big ),\quad \\
\rv(0)=\rv_{0},
\end{array}
\right .
\label{eq:Cuachyproblem}
\end{equation}
whose positive and continuous global solution is represented, thanks to the Leibnitz's formula,  by
$$
\int_{\rv_{0}}^{\rv(t)}\dfrac{ds}{\theta (s)}=\alpha t,\quad 0\le t\le \rT.
$$
When $\rv(t)$ satisfies the inequality \eqref{eq:maininequation}, we require a sharper refinement because an equivalence as  \eqref{eq:mainequation} and \eqref{eq:Cuachyproblem} does not hold in general. So, we introduce the positive and non decreasing  function $\disp \rV(t)=\max_{0\le \tau\le t} \rv(\tau)=\rv(\tau_{t})$, for some  $\tau_{t}\in [0,t]$. Next, we define
$$
\widehat{\rV}(t)= \rv_{0}+\alpha \int_{0}^{t}\theta \big (\rV(s)\big )ds > 0,\quad t>0,
$$
that satisfies $\widehat{\rV}(0)=\rv_{0}$, as well as
$$
\rV(t)=\rv(\tau_{t})\le \rv_{0}+\alpha \int_{0}^{\tau_{t}}\theta \big (\rv(s)\big )\le \rv_{0}+\alpha \int_{0}^{t}\theta \big (\rV(s)\big )=\widehat{\rV}(t) 
$$
and
$$
\left \{
\begin{array}{l}
\widehat{\rV}'(t)=\alpha \theta\big (\rV(t)\big )\le \alpha \theta\big (\widehat{\rV}(t)\big ),\\ [.15cm]
\widehat{\rV}(0)=\rv_{0}>0
\end{array}
\right .
$$
quite similar to \eqref{eq:Cuachyproblem}. Hence a kind of Leibnitz inequality
$$
\int^{\rv(t)}_{\rv_{0}}\dfrac{dr}{\theta(r)}\le \int^{\widehat{\rV}(t)}_{\rv_{0}}\dfrac{dr}{\vartheta(r)}=\int ^{t}_{0}\dfrac{\widehat{\rV}'(t)dt}{ \theta\big (\widehat{\rV}(t)\big )}\le \alpha t<+\infty
$$
holds, and then \eqref{eq:Leibnitzproperty} follows.
\par
\noindent On the other hand, when $\rv_{0}=0$, if we suppose $\rv(t)>0$ in some interval $t\in ]0 ,t_{1}]\subset [0,\rT]$ the above reasoning shows that $\widehat{\rV}(0)=0,~0<\rV(t)\le \widehat{\rV}(t)$ and
$ \widehat{\rV}'(t)\le \alpha \theta\big (\widehat{\rV}(t)\big )$, from which we deduce that $\widehat{\rV}(t)>0$ in $t\in ]0 ,t_{1}]$ and 
$$
\int^{\widehat{\rV}(t_{1})}_{0}\dfrac{dr}{\theta(r)}=\int ^{t_{1}}_{0}\dfrac{\widehat{\rV}'(t)dt}{\theta\big (\widehat{\rV}(t)\big )}\le \alpha t<+\infty,
$$
contrary to the condition \eqref{eq:Osgoodcondition}.$\fin$
\par
\medskip
Since $\theta$ is continuous and increasing we may introduce the increasing function
\begin{equation}
\Psi_{\rv_{0}}(\rv)\doteq \int^{\rv}_{\rv_{0}}\dfrac{ds}{\theta (s)},\quad \rv\ge \rv_{0},
\label{eq:Psifunction}
\end{equation} 
provided $\rv_{0}>0$. Then the inequality \eqref{eq:Leibnitzproperty} can be rewritten as 
\begin{equation}
\rv(t)\le \Psi ^{-1}_{\rv_{0}}\left (\alpha t\right ),\quad 0\le t\le \rT,
\label{eq:solutionintegralequation}
\end{equation} 
provided $\rv(0)=\Psi _{\rv_{0}}^{-1}(0)>0$. We emphasize that if $\rv_{0}=0$, the inequality \eqref{eq:solutionintegralequation} has not sense because the unique non negative function solving \eqref{eq:nullinequation} is the constant function $\rv(t)\equiv0$. It is easy to see that the function given by \eqref{eq:theta} verifies the Osgood assumption  \eqref{eq:Osgoodcondition}
$$
\int_{0^{+}}\dfrac{ds}{\theta (s)}=+\infty.
$$
\begin{rem}\em
Several generalizations of the above result will be presented in [\cite{diaz25abstract,diaz25deterministic}]. We note that in general no convex function $\theta(u)$ satisfies \eqref{eq:Osgoodcondition}. Certainly, if \eqref{eq:Osgoodcondition} holds the function $\theta$ is not integrable near 0. For instance the functions satisfying 
\begin{equation}
\dfrac{\theta(u)}{u}\le  \ln \dfrac{1}{u},\quad  \ln\dfrac{1}{u} \ln \stackrel{(n}{\cdots}\ln\dfrac{1}{u},\quad n\ge 0,\quad \hbox{near $u=0$},
\label{eq:mainexampleOsgood}
\end{equation}
provide other examples for which \eqref{eq:Osgoodcondition} holds near the origin. The conditions \eqref{eq:Osgoodcondition} and \eqref{eq:mainexampleOsgood} coincide with the classical Osgood's criterion (see \cite{osgood1898beweis}). On the other hand, we also note that the examples given by \eqref{eq:mainexampleOsgood} verify
$$
\theta (u)|\ln u|\le u |\ln u|\ln \dfrac{1}{u}\ln \stackrel{(n}{\cdots}\ln\dfrac{1}{u}=u\big (\ln u\big )^{n+2},\quad n\ge 0,\quad \hbox{near $\rv=0$},
$$
and thus they satisfy the so-called Dini condition
\begin{equation}
\lim _{u \searrow 0}\theta (u)|\ln u=0.
\label{eq:Dinicondition}
\end{equation}
\fineqnum
\end{rem}
We also emphasize that the inequality
\begin{equation}
|\theta(u)-\theta(v)|\le \theta(|u-v|),\quad u,v\ge 0.
\label{eq:modulusproperty}
\end{equation}
holds (see Remark \ref{rem:subadditiveproperty} below).

\begin{rem}\em In order to provide inequalities as \eqref{eq:modulusproperty} we may consider real functions $\Phi:[0,\infty[\rightarrow [0,\infty[$ satisfying $\Phi(0)\ge 0$ with  the subadditive property
\begin{equation}
\Phi(u)+\Phi(w)\ge \Phi(u+w)\quad \Leftrightarrow \quad \Phi(u+w)-\Phi(u)\le \Phi(w),\quad u,w\ge 0.
\label{eq:subadditive}
\end{equation}
holds. So that, let $u,v\ge 0$ such that $v\ge u$. We define $w=v-u\ge 0$ for which
$$
0\le \Phi(v)-\Phi(u)=\Phi(w+u)-\Phi(u)\le \Phi(w)=\Phi(v-u),
$$
whenever $\Phi $ is nondecreasing. By means of a similar reasoning, we conclude that
\begin{equation}
|\Phi (u)-\Phi (v)|\le \Phi (|u-v|),\quad u,v\ge 0,
\label{eq:continuitymodulus}
\end{equation}
whenever $\Phi $ is nondecreasing.$\fin$
\label{rem:continuitymodulus}
\end{rem}
\begin{rem}\em The sub-additive property \eqref{eq:subadditive} is satisfied by real concave functions $\Phi:[0,\infty[\rightarrow [0,\infty[$ satisfying $\Phi(0)\ge 0$. Indeed, it follows
$$
\Phi(\lambda z)=\Phi \big (\lambda z+(1-\lambda )0\big )\ge \lambda \Phi(z)+(1-\lambda )\Phi(0)\ge  \lambda \Phi(z),\quad z\ge 0,~0<\lambda<1.
$$ 
In particular, given $u,w>0$ by choosing $\lambda _{u,w}=\dfrac{u}{u+w} \in ]0,1[$ it follows
$$
\left \{
\begin{array}{l}
\Phi(u)=\Phi\big (\lambda _{u,w}(u+w)\big )\ge \lambda _{u,w}\Phi(u+w),\\
\Phi(w)=\Phi\big ((1-\lambda _{u,w})(u+w)\big )\ge (1-\lambda _{u,w})\Phi(u+w),
\end{array}
\right .
$$
whence one concludes the subadditive property
$$
\Phi(u)+\Phi(w)\ge \Phi(u+w).
$$
Arguing as in \cite[Lemma 4.3]{diaz2023large}, we also may obtain the subadditive property \eqref{eq:subadditive} by transfer without concavity settings. Indeed, assume
$$
q(u)+q(v)\ge q(u+v),\quad u,v\ge 0
$$
and
\begin{equation}
\dfrac{\Phi(u)}{q(u)}\quad \hbox{is non increasing}.
\label{eq:subadditivetransfer}
\end{equation}
Then
$$
\Phi(u)+\Phi(v)=q(u)\dfrac{\Phi(u)}{q(u)}+q(v)\dfrac{\Phi(v)}{q(v)}\ge (q(u)+q(v))\dfrac{\Phi(u+v)}{q(u+v)} \ge \Phi(u+v),
$$
thus, the  sub-additivity of the function $q(u)$ is transferred to the function $\Phi(u)$ provided \eqref{eq:subadditivetransfer}. In particular, any function $\Phi$ such that
$$
\dfrac{\Phi(u)}{u^{m}}\quad \hbox{is non increasing}
$$
for some $0<m\le 1$ is sub-additive and the inequality \eqref{eq:continuitymodulus}
holds whenever $\Phi $ is nondecreasing.$\fin$
\label{rem:subadditiveproperty}
\end{rem}
Finally,  we come back to the hybrid co-albedo profile (centered at the origin) of the Introduction
\begin{equation}
\beta(u)=\left\{ 
\begin{array}{ll}
\beta _{i}, & \quad \hbox{if $u<0$,} \\[0.1cm] 
\big (\beta _{w}-\beta_{i}\big )\dfrac{u\ln u}{\delta \ln \delta }+\beta_{i}, & \quad \hbox{if $0\le u\le \delta $}, \\[0.1cm] 
\beta _{w}, & \quad \hbox{if $u>\delta$},
\end{array}
\right.  
\label{eq:coalbedoprofile}
\end{equation}
with $0<\delta<1$ that governs the co-albedo function (see Figure 	\ref{fig:coalbedoprofile}).  We note that the profile verifies $\beta\in\cC(\RR)\cap \cC^{1}\big (\RR\setminus\{0,\delta\}\big )$, with
$$
\left\{
\begin{array}{l}
\beta'(0^{-})=0\quad \hbox{and \quad }\beta'(0^{+})=+\infty,\\
\beta'\big (\delta^{-}\big )>0\quad \hbox{and} \quad \beta'(\delta^{+})=0.
\end{array}
\right .
$$
We claim that the profile $\beta $ satisfies
\begin{equation}
|\beta(u)-\beta(v)|\le \theta(|u-v|),\quad u,v\in\RR.
\label{eq:moduloscontinuitybeta}
\end{equation}
(see \eqref{eq:modulusproperty}). Indeed, when $v\le 0$ one has
$$
\beta (u)-\beta (v)=\theta(u)-\theta(0)=\theta(u)\le  \theta(u-v) ,\quad v\le 0\le u\le \delta.
$$ 
Analogously, if $u\ge \delta$ one has
$$
\beta (u)-\beta (v)=\theta(\delta)-\theta(v)\le \theta(\delta-v)\le  \theta(u-v) ,\quad u\ge \delta\ge v\ge 0.
$$ 
Finally, the claim follows from
$$
\left \{
\begin{array}{l}
\beta (u)-\beta (v)=\theta(u)-\theta(v)\le \theta(u-v),\quad \delta \ge u\ge v\ge 0,\\ [.2cm]
\beta (u)-\beta (v)=\theta(\delta)-\theta(0)=\theta (\delta)\le \theta (u-v),\quad u\ge \delta \ge 0\ge v.
\end{array}
\right .
$$
Moreover, from \eqref{eq:Dinicondition} the function $\beta$ satisfies the Osgood's criterion \eqref{eq:Osgoodcondition} (see also Theorem \ref{theo:coalbedohybridfunction}). 

\section{On the cylindrical well-adapted to the Legendre diffusion Wiener noise}
\label{sec:preliminaries}
\par
By introducing the change of variable $\rX_{t}=\ru_{t}-10$,  $\xi=\ru_{0}-10$, the problem $(\rE_{\beta ,\varepsilon})$ corresponds to a choice of the general semilinear equation
\begin{equation}
\left \{
\begin{array}{l}
d \rX_{t}+\cA\rX_{t}dt=\rF_{t}\big (\rX_{t}\big )dt+\rB_{t}(\rX_{t})d \rW_{t},\\
\rX_{0}=\xi,
\label{eq:abstractCauchy}
\end{array}
\right .
\end{equation}
which will be treated on a separable Hilbert space $\rH$ where we are considering the measurable processes $\rX:\Omega_{\rT}\rightarrow \rH$, $\Omega_{\rT}=[0,\rT]\times \Omega$, posed in a probability space $(\Omega,\cF,\PP)$, equipped with a complete right continuous filtration $\{\cF_{t}\}_{t\in [0,\rT]}\subset \cF$. Moreover, as usual, we denote by  $\cP_{\rT}=\cB([0,\rT])\otimes \cF$ the predictable $\sigma$- fields on defined by
$$
\cP_{\rT}\doteq \sigma \left (\big \{]s,t]\times \rA_{s}:~0\le s<t\le \rT,~\rA_{s}\in \cF_{s}\big\}\cup \{\{0\}\times \rA_{0}: \rA_{0}\in \cF_{0}\}\right ).
$$
They may also be introduced by
$$
\cP_{\rT}\doteq \sigma \left ({\rY:~\Omega_{\rT}\rightarrow \RR:~\rY\hbox{ is left continuous and adapted  to }\cF_{t},~t\in [0,\rT]}\right ),
$$
more according to our purposes (see below).
\par
Next, we make precise the framework where \eqref{eq:abstractCauchy} will be formulated. Motivated by the Stochastic Partial Differential Equation of $(\rE_{\beta ,\varepsilon})$, we consider the differential operator $\rA:~\rD(\rA)\rightarrow \rH$ with
\begin{equation}
\left \{
\begin{array}{l}
\rD(\rA)=\big \{\rv\in \rH:~\rA \rv\in \rH\big\},\\ [.2cm]
\rA \rv(x)=-\dfrac{d}{dx}\left (\big (1-x^{2}\big )\dfrac{d}{dx}\rv(x)\right ),\quad x\in \rI=]-1,+1[,\quad \hbox{if $\rv\in \rD(\rA)$,}
\end{array}
\right .
\label{eq:Aoperator}
\end{equation}
defined in $\rH=\rL^{2}(\rI)$, equipped with its usual norm
$$
\n{\ru}=\left (\int_{\rI}|\ru(x)|^{2}dx\right )^{\frac{1}{2}} ,\quad \ru\in\rH
$$
and the scalar product 
$$
\pe{\ru}{\rv}=\int_{\rI}\ru(x)\rv(x)dx,\quad \ru,\rv\in\rH. 
$$
In \cite{D-Escorial} (see also \cite{Hetzer1990}) it was proved that $\rA$ is a maximal monotone operator densely defined, thus $\overline{\rD(\rA)}=\rH$.  More precisely,  according to \cite{D-Escorial} or \cite{GD-JID2022}, the domain of the operator is associated with a suitable energy space related to the Legendre diffusion operator. It is the weighted space
$$
\rV=\big \{\rw\in\rL^{2}(\rI):~\rw'\in \rL^{2}(\rI;\rho)\big \}
$$
where 
$$
\rL^{2}(\rI:\rho)=\left \{\rw:~\int _{\rI}\rho |\rw|^{2}dx<+\infty\right \},
$$
equipped with its norm
$$
\n{\rw}_{\rL^{2}(\rI:\rho)}^{2}=\int _{\rI}\rho |\rw|^{2}dx
$$
(we recall that $\rho(x)=1-x^{2}$). Here,  $\rH=\rL^{2}(\rI)$ is the so-called Hilbert pivot space.
Notice that $\rV$ is a separable Hilbert space related to the norm
$$
\n{\rw}_{\rV}=\n{\rw}_{\rL^{2}(\rI)}+\n{\rw'}_{\rL^{2}(\rI:\rho)}.
$$
Next, we introduce the abstract version of the diffusion operator by means of the functional operator $\cA:~\rV\rightarrow \rV'$ given by
\begin{equation}
\cA \ru\doteq -\dfrac{\partial }{\partial x}\left( \rho \dfrac{\partial }{\partial x}\ru\right),\quad \ru\in\rV.
\label{eq:leadingpart}
\end{equation}
Working with semigroup theory (\cite{brezis1973ope}), it is useful to define the above operator $\rA$ as the realization
$$
\left \{
\begin{array}{l}
\rD(\rA)=\big \{\rv\in \rH:~\cA \rv\in \rH\big\},\\ [.2cm]
\rA \rv(x)=\cA \rv\quad \hbox{if $\rv\in \rD(\rA)$.}
\end{array}
\right .
$$
Then, it was shown in \cite{D-Escorial} that the  the operator $\rA$ can be written as the subdifferential $\rA \rv =\partial \varphi (\rv)$
of the convex and lower semicontinuous functional, $\varphi :\rH\rightarrow \RR\cup \{+\infty \},$ given by 
\begin{equation}
\varphi (\rv)=\left\{ 
\begin{array}{ll}
\displaystyle\frac{1}{2}\int_{\rI}\rho (x)\left\vert \dfrac{\partial \rv}{\partial x}\right\vert ^{2}dx & \mbox{if }\bv\in \rV, \\ [.35cm]
+\infty  & \mbox{otherwise},%
\end{array}%
\right.   \label{eq:defphi}
\end{equation}
that $\rA$ is densely defined (see \cite[Proposition 1]{D-Escorial}) and that $\partial \varphi (\bv)$ generates a compact semigroup of contractions on $\rH$ (see \cite[Lemma 1]{D-Escorial}).
Then, for any $\ru_{0}\in\rH$ there exists a unique function $\ru\in\cC([0,\rT]:\rH)$, with the smoothing effect that $\ru(t)\in\rD(\rA)$ for $a.e.~ t\in (0,\rT] $, such that $\ru$ is the mild solution of the abstract problem   
\begin{equation}
\left\{ 
\begin{array}{l}
\dfrac{d \ru}{d t}(t)+\rA \ru(t)=0,\quad t>0,\\[0.2cm] 
\ru(0)=\ru_{0}\in\rH.
\end{array}
\right.  
\label{eq:deterabstractproblem}
\end{equation}
\begin{rem}\rm In fact, in \cite{D-Escorial} it was obtained the complementary regularity 
$$
t^{\frac{1}{2}}\dfrac{d \ru(t)}{dt}\in \rL^{2}\big (0,\rT:\HH\big )
$$ 
of the above mild solution. Moreover, if $\ru_{0}\in\rL^{p}(\rI),~1\le p\le \infty,$ then $\ru(t)\in\rL^{p}(\rI)$ for $a.e.~ t\in (0,\rT] $. In fact, if $\ru_{0}\in\rV\subset \rH$ one has $\dfrac{d \ru}{dt}\in \rL^{2}\big (0,\rT:\rH\big ).\fin$ 
\end{rem}

On the other hand, we also know that $\rH$  admits a Hilbertian basis given by the eigenvectors $\{\re_{n}\}_{n \ge 0}\subset\rD(\rA)$ of the operator $\rA$, defined through the  orthonormal Legendre polynomials of degree $n$, defined by the property
$$
\rA \rP_{n}=n(n-1)\rP_{n} \quad \hbox{on $\rH$}.
$$
So that, the constants $\mu_{n}=n(n+1)$ are the corresponding eigenvalues. Since we have
$$
\pe{\rP_{n}}{\rP_{m}}_{\rH}=\dfrac{2}{2n+1}\delta_{n,m},
$$
the normalized eigenvectors of $\rA$ are the functions
$$
\re_{n}(x)=\sqrt{\dfrac{2}{2n+1}}\rP_{n}(x),\quad -1\le x\le 1.\\ [.3cm]
$$
It also follows the Rodrigues formula that
$$
\rP_{n}(x)=\dfrac{1}{2^{n}n!}\dfrac{d^{n}}{dx ^{n}}\big (x^{2}-1\big )^{n}=\dfrac{1}{2^{n}}\sum_{k=0}^{n}\binom{n}{k}^{2}(x+1)^{n-k}(x-1),
$$
as well as the recurrence identity
$$
(n+1)\rP_{n+1}(x)+n\rP_{n-1}(x)=(2n+1)x\rP_{n}(x),
$$
(see, e.g., \cite{NorthLibro}).  Then, from the direct computation $\rP_{0}(x)\equiv 1$ and $\rP_{1}(x)=x$, we deduce that $\rP_{2}(x)=\dfrac{3x^{2}-1}{2},~\rP_{3}(x)=\dfrac{x(5x^{2}-3)}{2},\ldots $ and so on (see, e.g., \cite{nikiforov1991classical}). 
\begin{figure}[htp]
\begin{center}
\vspace*{.5cm}
\includegraphics[width=7cm]{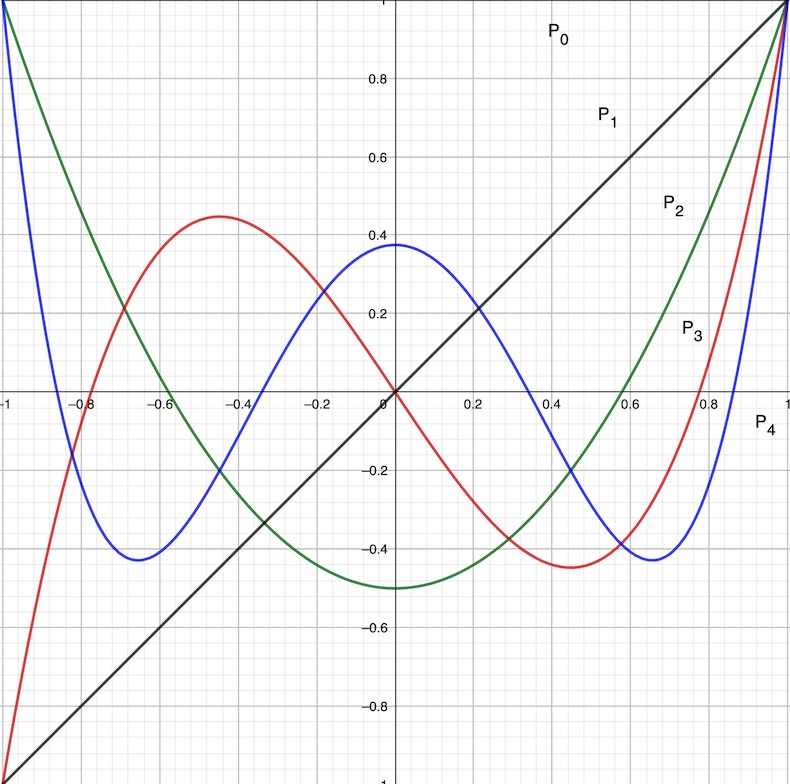}\\ 
\caption{Legendre polynomials}
\label{fig:Legendrepolynomials}
\end{center}
\end{figure}
\par
As mentioned in the Introduction, we emphasize that the first eigenvalue is the null value $\mu_{0}=0$. This implies that $\rA$ is not invertible. In order to avoid loss of invertibility, given $\mu >0$, we replace the differential operator $\rA$ by
$$
\rA_{\mu}\ru=\rA\ru+\mu\ru,\quad \ru\in\rD(\rA).
$$ 
Then, we get the same Hilbertian basis of $\rH$, given by the same eigenvectors $\{\re_{n}\}_{n \ge 0}\subset\rD(\rA)$ of the operator $\rA_{\mu}$
$$
\re_{n}(x)=\sqrt{\dfrac{2}{2n+1}}\rP_{n}(x),\quad -1\le x\le 1,
$$
but now, the corresponding eigenvalues are $ \widehat{\mu}_{n}=\mu_{n}+\mu>0$.
So that, for any $\rf\in\rH$ there exists a unique solution of $\rA_{\mu}\rv =\rf $ 
given by the representation
$$
\rv=\sum_{n\ge 0}\dfrac{\pe{\rf}{\re_{n}}}{\widehat{\mu}_{n}}\re_{n}.
$$
\par
\noindent
Denoting by $\{\rS^{\mu}_{t}\}_{t\ge 0}$ the semigroup generated by $\rA_{\mu}$,  for every $\ru_{0}\in \rH$, the mild solution of 
$$
\left\{ 
\begin{array}{l}
\dfrac{d \ru}{d t}(t)+\rA_{\mu} \ru(t)=0,\quad t>0,\\[0.2cm] 
\ru(0)=\ru_{0},
\end{array}
\right. 
$$
admits the representation
\begin{equation}
\ru(t)=\rS^{\mu}_{t}\ru_{0}=e^{-\mu t}\sum_{n\ge 0}\pe{\ru_{0}}{\re_{n}}_{\rH}\re^{-\mu_{n}t}\re_{n},\quad t\ge 0.
\label{eq:deterabstractproblemsolution}
\end{equation}
If we consider now the non-autonomous problem
\begin{equation}
\left\{ 
\begin{array}{l}
\dfrac{d \ru}{d t}(t)+\rA_{\mu} \ru(t)=\rf(t),\quad t>0,\\[0.2cm] 
\ru(0)=\ru_{0}\in\rH,
\end{array}
\right. 
\label{eq:deterabstractproblemf}
\end{equation}
when $\rf\in\rL^{1}\big (0,\rT;\rD(\rA)\big )$, we may solve \eqref{eq:deterabstractproblemf} via the generalized Duhamel formula (or constants variations formula) (see \cite{brezis1973ope})
$$
\ru(t)=\rS^{\widehat{\mu}}_{t}\ru_{0}+\int^{t}_{0}\rS^{\widehat{\mu}}_{t-s}\rf(s)ds,\quad 0\le s<\rT,
$$
thus
\begin{equation}
\ru(t)=\sum_{n\ge 0}\left (\pe{\ru_{0}}{e_{n}}_{\rH}e^{-\widehat{\mu}_{n}t}+e^{-\widehat{\mu}_{n}t}\int^{t}_{0}
\pe{\rf (s)}{\re_{n}}_{\rH}e^{\widehat{\mu} _{n}s}ds\right )\re_{n},\quad t\ge 0,
\label{eq:deterabstractproblemsolutionf}
\end{equation}
and $\ru$ is called the mild solution of \eqref{eq:deterabstractproblemf}.
\par
In order to make precise the stochastic framework, we introduce the operator $\cQ_{\mu}:\rH \stackrel{\rA_{\mu}^{-1}}{\rightarrow} \rD(\rA)\hookrightarrow\rH$, thus $\cQ_{\mu}\ru \in\rH$ is the solution of $\rA_{\mu}\cQ_{\mu}\ru =\ru $ with
$$
\cQ_{\mu} \ru=\sum_{n\ge 0}\dfrac{\pe{\ru}{\re_{n}}}{\widehat{\mu}_{n}}\ru_{n},\quad \ru\in\rH
$$
for which
$$
\n{\cQ_{\mu} \ru}^{2}=\sum_{n\ge 0}\dfrac{1}{\widehat{\mu}_{n}^{2}}\pe{\ru}{\re_{n}}^{2},\quad \ru \in\rH.
$$
Since 
$$
\sum_{n\ge 0}\dfrac{1}{\widehat{\mu}_{n}^{2}}<\sum_{n\ge 1}\dfrac{1}{\mu_{n}^{2}}< \sum_{n\ge 1}\dfrac{1}{\mu}_{n}=\sum_{n\ge 1}\left (\dfrac{1}{n}-\dfrac{1}{n+1}\right ) =1-\lim_{n\rightarrow \infty}\frac{1}{n+1}=1,
$$
the operator $\cQ_{\mu}$ is a bounded, positive and symmetric operator  with a finite trace
$$
{\rm Trace}~\cQ_{\mu}=\sum_{n\ge 1}\pe{\cQ_{\mu}\re_{n}}{\re_{n}}=\sum_{n\ge 0}\dfrac{1}{\widehat{\mu}_{n}}<\dfrac{1}{\mu}+\sum_{n\ge 1}\dfrac{1}{\mu}_{n}=\dfrac{1}{\mu}+1.
$$
\begin{rem}\rm
Since $\cQ_{\mu}\re_{n}=\dfrac{1}{\widehat{\mu}_{n}}\re_{n},~n\ge 0$ the  eigenvalues $\widehat{\lambda}_{n}$ of $\cQ_{\mu}$  coincide with the inverse values $\dfrac{1}{\widehat{\mu}_{n}}$ of those of the differential operator $\rA_{\mu}.\fin$ 
\end{rem}
\begin{rem}\rm
If we consider the nonnegative square root operator of $\cQ_{\mu}$,  the inequality 
$$
\n{\cQ_{\mu}^{\frac{1}{2}}\ru- \sum_{k=1}^{j}\pe{\cQ_{\mu}^{\frac{1}{2}}\ru}{\re_{k}}\re_{k}}^{2}\le\sum_{k=j+1}^{\infty}|\pe{\cQ_{\mu}^{\frac{1}{2}}\ru}{\re_{k}}|^{2}\le \n{\ru}^{2}\sum_{k=j+1}^{\infty}\n{\cQ_{\mu}^{\frac{1}{2}}\re_{k}}^{2}\le 
\n{\ru}^{2}\sum_{k=j+1}^{\infty}\pe{\cQ_{\mu}\re_{k}}{\re_{k}},
$$
shows that $\cQ_{\mu}^{\frac{1}{2}}$ is a compact operator because is a limit of finite rank operators, whence $\cQ_{\mu}=\cQ_{\mu}^{\frac{1}{2}}\cQ_{\mu}^{\frac{1}{2}}$ is also a compact operator. Moreover,
$$
\left \{
\begin{array}{l}
\disp \cQ _{\mu}\ru=\sum_{n\ge 0}\dfrac{1}{\widehat{\mu}_{n}}\pe{\ru}{\re_{n}}\re_{n},\\
\disp \cQ_{\mu}^{\frac{1}{2}}\ru = \sum_{n\ge 0}\dfrac{1}{\sqrt{\widehat{\mu}_{n}}}\pe{\ru}{\re_{n}}e _{n},
\end{array}
\quad \ru\in\rH.
\right .
$$
Finally, 
$$
\n{\cQ_{\mu}^{\frac{1}{2}}}_{\cL_{2}(\rH)}^{2}=\sum_{n\ge 0}\n{\cQ_{\mu}^{\frac{1}{2}}\re_{n}}^{2}=\sum_{n\ge 0}\dfrac{1}{\widehat{\mu}_{n}}=\n{\cQ_{\mu}}_{\cL_{1}^{+}(\rH)}=
\hbox{Trace }\cQ_{\mu}<+\infty,
$$
we have proved that $\cQ_{\mu}^{\frac{1}{2}}$ is a Hilbert-Schmidt operator on $\rH$.
\end{rem}
\begin{rem}\rm We note  the well-known representation of those Hilbert-Schmidt operators 
$$
\cQ_{\mu}^{\frac{1}{2}}\ru(x)=\int_{\rI}\left (\sum_{n\ge 0}\dfrac{1}{\sqrt{\widehat{\mu}_{n}}}\re_{n}(x)\re_{n} (y)\right )\ru(y)dy,\quad x\in\rI,\quad \ru\in\rH,
$$
whose kernel is
$$
k(x,y)\doteq \sum_{n\ge 0}\dfrac{1}{\sqrt{\widehat{\mu}_{n}}}\re_{m}(x)\re_{n} (y)\quad \hbox{with}\quad \int_{\rI\times\rI}|k(x,y)|^{2}dxdy\le \sum_{n\ge 0}\dfrac{1}{\widehat{\mu}_{n}}={\rm Trace}~\cQ_{\mu}.
$$
\fineq
\end{rem}
The operator $\cQ_{\mu}$ determines the so-called Cameron-Martin space $\rH _{\cQ_{\mu}}=\cQ_{\mu} ^{\frac{1}{2}}\rH$,   proper subspace dense in $\rH$, that
enables us  the relation 
$$
\ru=\cQ_{\mu}^{-\frac{1}{2}}\rv \quad \Leftrightarrow \quad \cQ_{\mu}^{\frac{1}{2}}\ru=\rv,
$$
thus the operator $\cQ_{\mu}^{-\frac{1}{2}}$ is defined on the Hilbert subspace $\rH_{\cQ_{\mu}}$ endowed with
$$
\n{\ru}_{\cQ_{\mu}}\doteq \n{\cQ_{\widehat{\mu}}^{-\frac{1}{2}}\ru}_{\rH},
$$
coming from the identity
$$
\pe{\ru}{\rv}_{\cQ_{\mu}}\doteq \pe{\cQ_{\mu}^{-\frac{1}{2}}\ru}{\cQ_{\mu}^{-\frac{1}{2}}\rv},\quad \ru,\rv\in \rH_{\cQ_{\mu}}.
$$
\begin{rem}\rm If there exists $\cQ_{\mu}^{-\frac{1}{2}}\ru=\rv\in\rH$ the property
$$
\cQ_{\mu}^{\frac{1}{2}}\rv=\ru\quad \Leftrightarrow \quad \dfrac{1}{\sqrt{\widehat{\mu}_{n}}}\pe{\rv}{\re_{n}}=\pe{\ru}{\re_{n}}
$$ 
requires  
$$
\n{\rv}^{2}=\sum_{n \ge 0}\widehat{\mu}_{n}\big |\pe{\ru}{\re_{n}}\big |^{2}<+\infty.
$$
Since $\ru\in\rH_{\cQ_{\mu}}$ implies $\ru=\cQ_{\mu}^{\frac{1}{2}}\rw$ for some $\rw\in\rH$ one has
$$
\pe{\ru}{\re_{n}}=\dfrac{1}{\sqrt{\widehat{\mu}}_{n}}\pe{\rw}{\re_{n}}.
$$
Then
$$
\sum_{n \ge 0}\sqrt{\widehat{\mu}}_{n}\big |\pe{\ru}{\re_{n}}\big |^{2}=\sum_{n \ge 0}\big |\pe{\rw}{\re_{n}}\big |^{2}=\n{\rw}^{2}.
$$
Hence, we have the representation
$$
\left \{
\begin{array}{l}
\disp \cQ_{\mu}^{-\frac{1}{2}}\ru=\sum_{n\ge 0}\sqrt{\widehat{\mu}_{n}}\pe{\ru}{\re_{n}}\re_{m}, \\[.2cm]
\disp \pe{\ru}{\rv}_{\cQ_{\mu}}=\sum_{n\ge 0}\widehat{\mu}_{n}\pe{\ru}{\re_{m}}\pe{\rv}{\re_{n}}=\pe{\cQ_{\mu}^{-\frac{1}{2}}\ru}{\cQ_{\mu}^{-\frac{1}{2}}\rv},\\ 
\end{array}
,\quad \ru \in\rH_{\cQ_{\mu}}.
\right .
$$
In particular, $\sqrt{\widehat{\mu}_{n}}\re_{n}\in\rH_{\cQ_{\mu}}$ and $\nl{\sqrt{\widehat{\mu}_{n}}\re_{n}}_{\cQ_{\mu}}=\n{\re_{n}}_{\rH}=1.\fin$
\end{rem}
\par
\medskip
After the above notations and commentaries we consider  a $\cQ_{\mu}$-cylindrical Wiener process on $\rH$ denoted by $\{\rW_{t}\}_{t\ge 0}$ satisfying 
\begin{description}
\item[i)] $\{\rW_{t}\}_{t\ge 0}$ has continuous trajectories and $\rW_{0}=0$,
\item[ii)] $\{\rW_{t}\}_{t\ge 0}$ has independent increments and
$$
\cL \big (\rW_{t}-\rW_{s}\big )=\cN\big (0,(t-s)\cQ_{\mu}\big )=\cL \big (\rW_{t-s}\big ),\quad t\ge s\ge 0
$$
\item[ii)] $\cL \big (\rW_{t}\big )=\cL \big (-\rW_{t}\big ),~t\ge 0$,
\end{description}
(see \cite{da2014stochastic}).
From the above definition, $\{\rW_{t}\}_{t\ge 0}$ is a Gaussian process on $\rH$ with
$$
\EE[\rW_{t}]=0\quad \hbox{ and }\quad {\rm Cov}~\rW_{t}=t\cQ_{\mu},\quad t\ge 0.
$$
Then, we may introduce the Brownian motions
\begin{equation}
\bB_{t}^{n}\doteq \sqrt{\widehat{\mu}}_{n}\pe{\rW_{t}}{\re_{n}},\quad t\ge 0.
\label{eq:BrownianQ}
\end{equation}
Since
$$
\begin{array}{ll}
\EE\big  [\bB_{t}^{n}\bB_{s}^{n'}\big ]&\hspace{-.2cm} \disp =\sqrt{\widehat{\mu}_{n}\widehat{\mu_{n'}}}\EE\big [\pe{\rW_{t}}{\re_{m}}\pe{\rW_{s}}{\re_{m'}}\big ]\\ [.3cm]
&\hspace{-.2cm} \disp =\sqrt{\widehat{\mu}_{n}\widehat{\mu_{n'}}}\left (\EE\big [\pe{\rW_{t}-\rW_{s}}{\re_{n}}\pe{\rW_{s}}{\re_{n'}}\big ]+
\EE\big [\pe{\rW_{s}}{\re_{n}}\pe{\rW_{s}}{\re_{n'}}\big ] \right )\\ [.3cm]
&\hspace{-.2cm} \disp =\sqrt{\widehat{\mu}_{n}\widehat{\mu_{n'}}}s\pe{\cQ \re_{n}}{\re_{n'}}=s\delta_{n,n'},\quad 0\le s\le t,
\end{array}
$$
the processes $\big \{\bB_{t}^{m}\big \}_{t\ge 0}$ are mutually independent. Moreover $\EE\big [|\bB^{k}_{t}|^{2}\big ]=t$, uniformly on $k$, implies
$$
\EE \left [\nll{\sum_{k=j}^{n}\dfrac{1}{\sqrt{\widehat{\mu}_{n}}}\bB_{t}^{k}\re_{k}}^{2}\right ]=t\sum_{k=j}^{n}\pe{\cQ_{\mu}\re_{k}}{\re_{k}}=t\sum_{k=j}^{n}\dfrac{1}{\sqrt{\mu_{k}}}.
$$
Thus, the property ${\rm Trace}~\cQ_{\mu}<\infty $ enables us to admit the representation on $\rH$
\begin{equation}
\rW_{t}=\sum_{n\ge 0}\dfrac{1}{\sqrt{\widehat{\mu}}_{n}}\bB_{t}^{n}\re_{n},\quad t\ge 0,
\label{eq:WienerQ}
\end{equation}
for which the Wiener isometry
\begin{equation}
\EE \big [\n{\rW_{t}}^{2}\big ]=
\EE \left [\sum_{n\ge 0}\dfrac{1}{\widehat{\mu}_{n}}\n{\bB_{t}^{n}\re_{n}}^{2}\right ]=t\sum_{n\ge 1}\dfrac{1}{\widehat{\mu}_{n}}=t{\rm Trace}~\cQ_{\mu}
\label{eq:Wienerisometry1}
\end{equation}
holds. 
\begin{rem}\rm In fact, from the maximality martingale inequality (see \cite{da2014stochastic})
$$
\PP \left(\sup_{t\in [0,\rT]}\sum_{k=j}^{n}\dfrac{1}{\sqrt{\widehat{\mu}_{k}}}\n{\bB_{t}^{k}\re_{k}}>r\right)
\le \dfrac{4}{r} \EE\left [\nl{\sum_{k=j}^{n}\dfrac{1}{\sqrt{\widehat{\mu}_{k}}}\bB_{\rT}^{j}\re_{k}}^{2}\right ]\le \dfrac{4\rT}{r}\sum_{k=j}^{n}\dfrac{1}{\sqrt{\widehat{\mu}_{k}}},
$$
one proves that  the serie \eqref{eq:WienerQ} is uniformly convergent on  $[0,\rT]~\PP$-$a.s$ (see \cite[Theorem 4.3]{da2014stochastic}).$\fin$
\end{rem}
\begin{rem}\rm The isometry \eqref{eq:Wienerisometry1}, as well as
$$
\EE \big  [\pe{\rW_{t}}{\ru}\pe{\rW_{s}}{\rv}\big ]=\EE \left [\sum_{n\ge 0}\dfrac{1}{\widehat{\mu}_{n}}\bB^{n}_{t}\bB^{n}_{s}\pe{\re_{n}}{\ru}\pe{\re_{n}}{\rv}\right ]=t\wedge s\pe{\cQ_{\mu}\ru}{\rv},
$$
show the covariance property of the operator $\cQ_{\mu}.\fin$
\label{rem:generalisometryW}
\end{rem}
\par
On the other hand, we denote by $\cL_{\cQ_{\mu}}=\cL_{2}(\rH_{\cQ_{\mu}},\rH)$ the space of  the Hilbert-Schmidt operators $\cD:~\rH_{\cQ_{\mu}}\rightarrow \rH$ equipped with the norm
$$
\n{\cD}_{\cL_{\cQ_{\mu}}}^{2}=\sum_{n\ge 0}\n{\cD\cQ_{\mu}^{\frac{1}{2}}\re_{n}}^{2}=
\sum_{n\ge 0}\n{\cD\cQ_{\mu}^{\frac{1}{2}}\cQ_{\mu}^{\frac{1}{2}}\cD^{*}\re_{n}}^{2}.
$$
In fact, one has  
$$
\n{\cD}_{\cL_{\cQ_{\mu}}}^{2}=\sum_{n\ge 0}\n{\cD\cQ_{\mu}^{\frac{1}{2}}\re_{n}}^{2}={\rm Trace}~\cD\cQ_{\mu}\cD^{*}=\sum_{n\ge 0}\dfrac{1}{\widehat{\mu}_{n}}\n{\cD\re_{n}}^{2}.
$$
So that, one defines the $\cQ_{\mu}$-predictable process as the process $\rB:\Omega_{\rT}\rightarrow \cL_{\cQ_{\mu}},~\rB _{t}(\cdot)(\omega)\in \cL_{\cQ_{\mu}}$ for which
$$
\nn{\rB}_{\cP_{\rT}}\doteq \left (\EE \left [\int^{\rT}_{0}{\rm Trace}~ \rB_{t}\cQ_{\mu}\rB^{*}_{t}dt\right ]\right )^{\frac{1}{2}}= \left (\EE \left [\int^{\rT}_{0}\n{\rB_{t}}_{\cL_{\cQ_{\mu}}}^{2}dt\right ]\right )^{\frac{1}{2}}<\infty.
$$ 
\par
Among the $\cQ_{\mu}$-predictable processes, we focus on the stochastic integral 
$$
\rB\cdot \rW\doteq \int^{\rT}_{0}\rB_{s}d\rW_{t},
$$  
where $\{\rB_{t}\}_{t\ge 0}\in \cP_{\rT}$ is a $\{\cF_{t}\}_{t\ge 0}$ adapted process $\cL_{\cQ_{\mu}}$ valued posed in $\rH_{\cQ_{\mu}}$ (see \cite{chow2007stochastic,da2014stochastic,liu2015stochastic} for definition) for which one introduces the processes $\{(\rB\cdot\rW)_{t}\}_{t\in [0,\rT]}$ given by
$$
(\rB\cdot\rW)_{t}\doteq \1_{[0,t]}\rB\cdot \rW
$$ 
(see \cite{chow2007stochastic,liu2015stochastic} again). Thus $(\rB\cdot \rW)_{t}(\omega)\in \cL_{\cQ_{\mu}},~\omega\in\Omega,$ satisfies the Ito isometry  
\begin{equation}
\EE\big [\n{(\rB\cdot\rW)_{\rT}}^{2}\big ]=\left (\EE \left [\int^{\rT}_{0}{\rm Trace }~ \rB_{t}\cQ_{\mu}\rB^{*}_{t}dt\right ]\right )^{\frac{1}{2}}= \left (\EE \left [\int^{\rT}_{0}\n{\rB_{t}}^{2}_{\cL_{\cQ_{\mu}}}dt\right ]\right )^{\frac{1}{2}}=\nn{\rB}_{\cP_{\rT}}^{2}<\infty.
\label{eq:isometryIto}
\end{equation}

\par
Relative to the other terms of the general semilinear problem
\begin{equation}
\left \{
\begin{array}{l}
d \rX_{t}+\rA_{\mu}\rX_{t}dt=\rF_{t}(\rX_{t})dt+\rB_{t}(\rX_{t})d \rW_{t},\\
\rX_{0}=\xi,
\end{array}
\right .
\label{eq:abstractCauchymu}
\end{equation}
on the Hilbert space $\rH$ (see \eqref{eq:abstractCauchy}), we will require that $\rF_{t}(\ru)(\omega)\in\rH,~(t,\omega,\ru) \in \Omega_{\rT}\times \rH$ for which  the Bochner integral
$$
\int^{\rT}_{0}\rF_{s}ds 
$$
posed in $\rH$ (see \cite{liu2015stochastic}) is well defined. Finally, once again, we recall that  $\rA_{\mu} $ is the infinitesimal generator of a strongly linear semigroup of contractions $\{\rS^{\mu}_{t}=e^{-\mu t}\rS_{t}\}_{t\ge 0}$ in $\overline{\rD(\rA_{\mu})}=\overline{\rD(\rA)}=\rH$.
\par
From the constants variation formula we get that a solution of \eqref{eq:abstractCauchymu} is a  measurable process $\rX$ from $(\Omega_{\rT},\cP_{\rT})$ into $(\rH,\cB_{\rH})$  satisfying for arbitrary $t\in [0,\rT]$
$$
\rX_{t}=\rS^{\mu}_{t}\xi +\int_{0}^{t}\rS^{\mu}_{t-s}\big (\rF_{s}(\rX_{s})\big) ds+\int_{0}^{t}\rS^{\mu}_{t-s} \rB_{s}(\rX_{s})d\rW_{s},\quad \PP~a.e.
$$
where $\xi$ is a  $\rH$-valued random variable. We will get later that the following property holds
$$
\PP \left (\int_{0}^{t}\bigg(\nl{\rS^{\mu}_{t-s}\rF_{s}(\rX_{s})}_{\rH}+\nl{\rS^{\mu}_{t-s} \rB_{s}(\rX_{s}) }^{2}_{\cL_{\cQ}}\bigg )ds<+\infty\right )=1,
$$
in order to the above equality is well defined. We call such a type of function $\rX$ a {\em mild solution}.  
\par
So, we will solve the semilinear integral equation by proving the existence of a process $\rX$ which is a fixed point of the operator $\cG$
$$
\rX=\cG^{\xi}\rX,
$$
where
$$
\big (\cG^{\xi} \rX\big )_{t}\doteq \rS^{\mu}_{t}\xi +\int_{0}^{t}\rS^{\mu}_{t-s} \rF_{s}\big (\rX_{s}\big )ds+\int_{0}^{t}\rS^{\mu}_{t-s}\rB_{s} \big (\rX_{s}\big )d\rW_{s}.
$$
Certainly, any such fixed point $\rX$ solves the semilinear integral equation. Among others possibilities, a way to prove the existence of a fixed point $\rX$ is to find some suitable topology in which  the relative Picard type of successive approximations
$$
\big (\rX_{0}\big )_{t}=\rS^{\mu}_{t}\xi  ,\qquad \rX_{n+1}=\cG^{\xi} \rX_{n},~n\ge 0,
$$
converge to $\rX$. This approximation problem will be used in the next Section \ref{sec:climatemodels}.

\begin{rem}[Doss-Sussman type transformations]\rm \label{Rmchangevariabl} There are some other ways to study semilinear stochastic equations.  This is the case of some useful transformations as the one due to H. Doss and H.J. Sussman \cite{doss1977liens,sussmann1978gap} (see also the pioneering transformation made in \cite{bensoussan1972equations}). In particular, in the multiplicative noise case, we may consider the transformation $\rX_{t}=\Psi(\rY_{t},\rW_{t}) $, for which Ito' Rule gives
\begin{equation}
\begin{array}{ll}
	d \rX_{t}-\rB (\rX_{t})d \rW_{t}&\hspace*{-.2cm}= \big [\Psi_{w}(\rY_{t},\rW_{t})-\rB \big (\Psi(\rY_{t},\rW_{t}\big )\big ]d\rW_{t}\\
	&+\Psi_{y}(\rY_{t},\rW_{t})d \rY_{t}+\dfrac{1}{2}\Psi_{ww}(\rY_{t}, \rW_{t})\sigma(\rW_{t})\sigma^{*}(\rW_{t})dt,
\end{array}
\end{equation}
whence
\begin{equation}
d \rX_{t}-\rB_{t} (\rX_{t})d \rW_{t}=\Psi_{y}(\rY_{t},\rW_{t})d \rY_{t}+\dfrac{1}{2}\Psi_{ww}(\rY_{t},\rW_{t})\sigma(\rW_{t})\sigma^{*}(\rW_{t})dt, 
\label{DossSussmanconversion}
\end{equation}
provided
\begin{equation}
\Psi_{w}(y,w)=\rB\big (\Psi(y,w)\big).
\label{DossSussman}
\end{equation}
Thus, by means of \eqref{DossSussman} the Doss-Sussman transfomation $\rX_{t}=\Psi(\rY_{t},\rW_{t})$ converts an stochastic integral in a random integral (see \eqref{DossSussmanconversion}). In particular, 
\begin{equation}
d \rX_{t}=\big (\cA \rX_{t}+\rF_{t}(\rX_{t})+\rf_{t}\big )dt+\rB(\rX_{t})d \rW_{t},\\
\label{eq:abstract1}
\end{equation}
becomes the random differential equation (which can be considered as a kind of deterministic nonlinear equation)
$$
\dfrac{d\rY_{t}}{dt}=\dfrac{1}{\Psi_{y}(\rY_{t},\rW_{t})}\left [\cA \Psi(\rY_{t},\rW_{t})+\rF_{t}\big (\Psi(\rY_{t},\rW_{t})\big )+\rf_{t}+\dfrac{1}{2}\Psi_{ww}(\rY_{t},\rW_{t})\sigma (\rW_{t})\sigma^{*}(\rW_{t})\right ].
$$
The study of \eqref{eq:abstract1} by the Doss-Sussman transformation is very tedious whenever  $\rB(\rX_{t})$ is a general nonlinear diffusion term.
\par
\noindent
Nevertheless, as mentioned in the Introduction, in the multiplicative linear case $\rB_{t}(\rX_{t})=a \rX_{t}$, we may take $\Psi(y,w)=y\Phi(w)$. Then the Doss-Sussman equation \eqref{DossSussman} becomes
$$
\Phi'(w)=a\Phi(w)\quad \Rightarrow \quad \Phi(w)=e^{a w}
$$
and we obtain the random differential equation
$$
\dfrac{d\rY_{t}}{dt}=e^{-a\rW_{t}}\cA e^{a\rW_{t}}\rY_{t}+e^{-a\rW_{t}}\big (\rF_{t}\big (e^{a\rW_{t}}\rY_{t}\big )+\rf_{t}\big )+\dfrac{1}{2}\rY_{t}\sigma (\rW_{t})\sigma^{*}(\rW_{t}).
$$
On the other hand, in the additive noise case $\rB_{t}(\rX_{t})\equiv \rB_{t}$ this type of Doss-Sussman transformation is quite simpler (see e.g., \cite{GD-JID2022}) 
$$
d \rX_{t}-\rB_{t}d \rW_{t}=d\rY_{t}\quad \Rightarrow \quad \rX_{t}=\rY_{t}+(\rB\cdot \rW)_{t},
$$
and then
$$
d\rY_{t}=d \rX_{t}-\rB d \rW_{t}=\big (\cA \rX_{t} +\rF_{t}(\rX_{t}\big )dt,
$$
which becomes the random differential equation
$$
\dfrac{d \rY_{t}}{dt}=\cA \rY_{t}+\rF_{t} (\rY_{t}+(\rB\cdot \rW)_{t})+\rf_{t}+\cA (\rB\cdot \rW)_{t}.
$$
\fineq
\end{rem}

\section{An application of the successive approximation method to the stochastic energy balance  model with a non-Lipschitz co-albedo}
\label{sec:climatemodels}
\par
In what follows, we assume the hypotheses (\bH$_{g}$) and (\bH$_{s}$), with $\beta $ given by (\ref%
{eq:coalbedofunction}), as in the Introduction. We will apply the stochastic framework of the Section \ref{sec:preliminaries}. In particular, we will use 
the notations and comments made there when we look at problem $(\rE
_{\beta ,\varepsilon})$ as a special case of the abstract stochastic equation 
\begin{equation}
d\ru_{t}+\big (\rA \ru_{t}+\rR_{e}\big )dt=\rR_{a}\big (dt+\varepsilon d\rW_{t}\big ),\quad t>0, 
\label{eq:SDEabs}
\end{equation}
prescribing the initial datum $\ru_{0}\in\rH$ (see \eqref{eq:abstractCauchy} where $\{\rW_{t}\}_{t \ge 0}$ was given in  \eqref{eq:WienerQ}).
\par
As it was pointed out, different kind of notions of solutions are possible, and then it is crucial to formulate correctly the assumptions on the data. At least formally, the problem  \eqref{eq:SDEabs} is equivalent to the integral identity 
\begin{equation}
\ru_{t}=\ru_{0}+\int_{0}^{t}(-\rA \ru_{s}+\big (\rR_{e}-\rR_{a}\big )ds+\varepsilon\int_{0}^{t}\rR_{a}d\rW_{s},\quad \forall t>0,
\label{eq:strongsolution} 
\end{equation}
where $\{\ru_{t}\}_{t\geq 0}$  must be an adapted random process to the filtration satisfying, in some sense, the integral representation \eqref{eq:strongsolution}. 
\par
\noindent
From Section \ref{sec:preliminaries} we recall that the operator $\rA:\rD(\rA)\rightarrow \rH$ defined in \eqref{eq:Aoperator} generates a semi-group $\{\rS_{t}\}_{t\ge 0}$ of contractions on $\rH$. In fact, since $\rA$ is not invertible, given $\mu>0$, we replace it by the operator $\rA_{\mu}\ru=\rA\ru+\mu\ru,$ and consider the equation
\begin{equation}
d\ru_{t}+\big (\rA_{\mu} \ru_{t}-\mu\rX_{t}+\rR_{e}\big )dt=\rR_{a}\big (dt+\varepsilon d\rW_{t}\big ),\quad t>0. 
\label{eq:SDEabsmu}
\end{equation}
On the other hand, as in the Budyko proposal, we assume that the Earth's radiation is of the type
$$
\rR_{e}(x,t,u)=g\big (u(x,t)\big )
$$
where $g:\RR\rightarrow \RR$ is a continuous increasing (see (\bH$_{g}$)) (see the Introduction).  The above change in the diffusion operator implies that now we include in \eqref{eq:SDEabs} an artificial term 
$$
\widehat{\rR}_{e}=g\big (\ru_{t}\big )-\mu\ru_{t}.
$$
The absorbed radiation energy is given by
$$
\rR_{a}(x,t,u)=\rQ\rS(x)\beta _{-10}\big (u(x,t)\big ),
$$
under the assumption $\big (\bH_{s}\big )$ and involving the hybrid co-albedo function  $\beta_{-10}(u)=\beta (u-10)$, where the function $\beta$ is the coalbedo profile given in \eqref{eq:coalbedoprofile}. 
\par
As before, we simplify the exposition by making the change of unknown $\rX_{t}=\ru_{t}-10$.  
Then, \eqref{eq:SDEabs} becomes
\begin{equation}
\left\{ 
\begin{array}{l}
d\rX_{t}+\rA_{\mu} \rX_{t}dt=\rF(\rX_{t})dt +\rB(\rX_{t})d\rW_{t},\quad t>0, \\
[0.1cm]
\rX_{0}=\ru_{0}-10\in \rH,
\end{array}
\right.   
\label{eq:SDEabsFG}
\end{equation}
where
\begin{equation}
\left \{
\begin{array}{l}
\rF(\rX_{t})=\mu\big (\rX_{t}+10) -g\big (\rX_{t}+10\big )+\rQ \rS \theta(\rX_{t}),\\ [.1cm]
\rB(\rX_{t})=\varepsilon\rQ\rS \theta(\rX_{t})\in\cL_{\cQ_{\mu}}^{2},
\end{array}
\right .
\label{eq:data}
\end{equation}
(see Section \ref{sec:preliminaries}).
\begin{rem}\rm For the condition $\rB(\rX)\in\cL_{\cQ_{\mu}}^{2}$ we mean that it is 
given by
$$
\big (\rB(\rX)u)(x)\doteq \varepsilon\rQ\rS(x) \theta(\rX)u(x),\quad x\in\rI,\quad \ru\in\rH.
$$
Since $\rS\in\rL^{\infty}(\rI)$  one has that $\rB(\rX)\ru\in\rH$ for $\ru\in\rH$. Moreover
\begin{equation}
\n{\rB(\rX)}_{\cL^{2}_{\cQ_{\mu}}}^{2}\le \varepsilon^{2}\rQ^{2} \big (\theta(\rX)\big )^{2}\n{\rS}_{\rL^{\infty}(\rI)}^{2}{\rm Trace }~\cQ_{\mu}.
\label{eq:boundbB}
\end{equation}
\fineqnum
\end{rem}
From the constants variation formula we have that any solution of \eqref{eq:SDEabsFG} must be a  measurable process $\rX$ from $(\Omega_{\rT},\cP_{\rT})$ into $(\rH,\cB_{\rH})$  satisfying, for arbitrary $t\in [0,\rT]$,
\begin{equation}
\rX_{t}=e^{-\mu t}\rS_{t}\rX_{0}+\int^{t}_{0}e^{-\mu (t-s)}\rS_{t-s}\rF\big(\rX_{s}\big )ds+\int^{t}_{0}e^{-\mu (t-s)}\rS_{t-s}\rB(\rX_{s})d\rW_{s},\quad t\ge 0.
\label{eq:mildsolutionmuFB}
\end{equation}
Notice that here $\disp \rM\doteq \sup_{0\le t\le \rT}e^{-\mu t}\n{\rS_{t}}_{\rL(\HH)}\le 1$.  We will also check that the following property holds
\begin{equation}
\PP \left (\int_{0}^{t}\bigg(\nl{e^{-\mu (t-s)}\rS_{t-s}\rF(\rX_{s})}_{\rH}+\nl{e^{-\mu (t-s)}\rS_{t-s} \rB(\rX_{s}) }^{2}_{\cL^{2}_{\cQ_{\mu}}}\bigg )ds<+\infty\right )=1.
\label{eq:semilinearintegralequation2}
\end{equation}
When \eqref{eq:mildsolutionmuFB} and \eqref{eq:semilinearintegralequation2} hold, we say that the process $\rX$ is a mild solution.  
\par
Due to the presence of nonlinear terms in \eqref{eq:SDEabsFG}, our goal is to prove the existence of processes $\rX^{\ru_{0}}$ which are fixed points of $\cG^{\ru_{0}}$,
$$
\rX^{\ru_{0}}=\cG^{\ru_{0}}\rX^{\ru_{0}},
$$
for the operator
\begin{equation}
\big (\cG^{\ru_{0}} \rX\big )_{t}\doteq e^{-\mu t}\rS_{t}\rX_{0}+\int^{t}_{0}e^{-\mu (t-s)}\rS_{t-s}\rF\big(\rX_{s}\big )ds+\int^{t}_{0}e^{-\mu (t-s)}\rS_{t-s}\rB(\rX_{s})d\rW_{s}.
\label{eq:functionaloperatorG}
\end{equation}
Among other possible arguments, a way to prove the existence of a fixed point $\rX^{\ru_{0}}$ is to find some suitable topology in which the relative Picard type of successive approximations
\begin{equation}
\big (\rX_{0}\big )_{t}= e^{-\mu t}\rS_{t}(\ru_{0}-10) ,\qquad \rX_{n+1}=\cG^{\ru_{0}} \rX_{n},~n\ge 0
\label{eq:Picardapproximations}
\end{equation}
converge to $\rX^{\ru_{0}}$.
Some reasons (see \eqref{eq:estimatekey2} below) advise introducing the Banach space (where we will solve the fixed-point problem) given by 
\begin{equation}
\BB_{\rT}\doteq \left \{\rX\in \cP_{\rT}:~\displaystyle\EE\left [\sup_{0\le s\le \rT}\n{\rX_{s}}_{\rH}^{2}\right ]<\infty \right \}
\label{eq:BTspace}
\end{equation}
endowed with the norm
$$
\n{\rX}_{\BB_{\rT}}\doteq \left (\EE\left [\sup_{0\le s\le \rT}\n{\rX_{s}}_{\rH}^{2}\right ]\right )^{\frac{1}{2}}.
$$
Certainly, we may then consider the subspaces  $\BB_{t}\subset \BB_{\rT}$ equipped with
$$
\n{\rX}_{\BB_{t}}\doteq \left (\EE\left [\sup_{0\le s\le t}\n{\rX_{s}}_{\rH}^{2}\right ]\right )^{\frac{1}{2}}\le
\n{\rX}_{\BB_{\rT}},
$$ 
for each $t\in [0,\rT]$.

Under Lipschitz assumptions on the nonlinear terms $\rF$ and $\rB$ the applicability of the classical successive approximation method is well known in the literature (see, e.g., \cite{da2014stochastic} and \cite{yamada1981successive}). When Lipchitz assumptions do not hold, other assumptions are needed. We send to \cite{barbu2002approximations,taniguchi1992successive} where the adaptation of the method under non Lipschitz conditions is considered.

A key stone in our reasoning is based on  an extension of the maximal sub-martingale inequality (see  \cite[Lemma 7.2]{da2014stochastic})  that we apply to the stochastic convolution  term
$$
\big (\rW_{\rA}^{\rB(\rX)}\big )_{t}\doteq \int_{0}^{t}\rS_{t-s}\rB_{s}\big (\rX_{s}\big )d\rW_{s}.
$$
More precisely, we will use a suitable simple consequence of \cite[Proposition 7.3]{da2014stochastic} (see also \cite[Theorem 1]{tubaro1984estimate}).
\begin{prop}[\cite{diaz25abstract}]. There exists a positive constant $\rc_{\rT}$ such that 
\begin{equation}
\EE \left [\sup_{0\le t\le \rT}\nll{\int_{0}^{t}\rS_{t-s}\rB\big (\rX_{s}\big )d\rW_{s}}^{2}_{\rH}\right ]\le \rc_{\rT}\int^{\rT}_{0}\EE \big [\nl{\rB\big (\rX_{s}\big )}^{2}_{\cL_{\cQ_{\mu}}}\big ]ds.
\label{eq:estimatekey2}
\end{equation}
holds. $\fin$
\label{prop:key2}
\end{prop}

\begin{rem}\rm The general estimate 
\begin{equation}
\EE \left [\sup_{0\le t\le \rT}\nll{\int_{0}^{t}\rS_{t-s}\rB\big (\rX_{s}\big )d\rW_{s}}^{p}_{\rH}\right ]\le \rc_{\rT}\EE \left [\int^{\rT}_{0}\nl{\rB\big (\rX_{s}\big )}^{p}_{\cL_{\cQ_{\mu}}}ds\right ],
\label{eq:estimatekeyp}
\end{equation}
only holds for some power $p>2$ (see \cite[Proposition 7.3]{da2014stochastic}).  However, when $\{\rS_{t}\}_{t\ge 0}$ is a semigroup of contractions, the case $p=2$ 
\begin{equation}
\EE \left [\sup_{0\le t\le \rT}\nll{\int_{0}^{t}\rS_{t-s}\rB\big (\rX_{s}\big )d\rW_{s}}^{2}_{\rH}\right ]\le \rc_{\rT}\EE \left [\int^{\rT}_{0}\nl{\rB\big (\rX_{s}\big )}^{2}_{\cL_{\cQ_{\mu}}}ds\right ]
\label{eq:estimatekey22}
\end{equation}
is valid (see \cite[Theorem 1]{tubaro1984estimate}).  In fact, arguing as  in \cite[Theorem 4.7]{da2014stochastic} we also arrive to \eqref{eq:estimatekey2}$\fin$
\end{rem}

\begin{figure}[htp]
\begin{center}
\includegraphics[width=3.75cm]{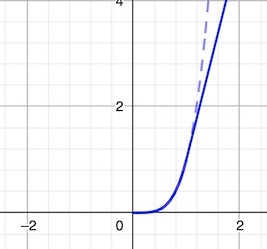}\\ 
\caption{Lipschitz truncation}
\label{fig:Lipschitztruncation}
\end{center}
\end{figure}

The approximation \eqref{eq:Picardapproximations} will be studied, in a more general framework, in \cite{diaz25abstract}. As in \cite{barbu2002approximations}, this idea was motivated by  \cite[Theorem 7.2]{da2014stochastic}. 
\begin{rem}\rm  Since our datum $u_{0}$ will be assumed to be a bounded function, it is natural to search for bounded solutions of the deterministic associated problem
$$
(\rE_{\beta ,0})
\left\{ 
\begin{array}{l}
\dfrac{\partial u}{\partial t}(x,t)-\dfrac{\partial }{\partial x}\left (\big (
1-x^{2})\dfrac{\partial u}{\partial x}\right  )(x,t)+g\big (u(x,t)\big )= 
\rQ\rS(x)\beta _{-10}\big (u(x,t)\big ),\\[0.2cm] 
u(x,0)=u_{0}(x),
\end{array}
\right. 
$$
when $g(u)$ is assumed to be locally Lipschitz continuous and increasing. In short, if we try to build a constant supersolution    
$$
\overline{u}(x,t)=\rK,
$$
then, we get that it is enough to assume
$$
\rK\ge \max \left \{ \n{u_{0}}_{\rL^{\infty}(\rI)}, g^{-1}\left (\rQ\n{\rS}_{\rL^{\infty}(\rI)}\beta_{w}\right )\right \}.
$$
(see \eqref{eq:coalbedoprofile}). This explains that we can assume that $g(u)$ is globally Lipschitz continuous, since we may replace it by the truncated function
\begin{equation} 
\widehat{g}(u)=
\left \{
\begin{array}{ll}
	g(-\rK)+g'(-\rK)(u-\rM)&\quad \hbox{if $u<\rK$},\\ [.1cm]
	g(u) &\quad \hbox{if $-\rK\le u\le \rK$},\\ [.1cm]
	g(\rK)+g'(\rK)(u+\rK)&\quad \hbox{if $\rK<u$},
\end{array}
\right .
\label{eq:truncationg}
\end{equation}	
(see Figure \ref{fig:Lipschitztruncation}).$\fin$ 
\label{rem:boundedness}
\end{rem}
Here we will apply the successive approximation method to the case of the hybrid co-albedo function \eqref{eq:theta}. By using Theorem \ref{theo:coalbedohybridfunction} and Lemma \ref{lemma:estimatekey2refined} below we will prove the following result:

\begin{prop} Let $\ru_{0},~\widehat{\ru}_{0}\in \rH$. Then the operator $\cG^{\ru_{0}}:~\BB_{\rT}\rightarrow \BB_{\rT}$ given by \eqref{eq:functionaloperatorG} is well defined and continuous. Moreover, one has the estimate
\begin{equation}
\disp \nl{\cG^{\ru_{0}}\rX-\cG^{\widehat{\ru}_{0}}\rY}_{\BB_{t}}^{2}\le
4\rM^{2}\n{\ru_{0}-\widehat{\ru}_{0}}_{\rH}^{2}+\widehat{\rC}_{\rT} \int^{\rT}_{0}\theta_{\rF,\rB}\big (\nl{\rX-\rY}_{\BB_{s}}^{2}\big )ds ,\quad t\in [0,\rT].
\label{eq:Gcontinuous}
\end{equation}
Here $\disp \rM\doteq \sup_{0\le s\le \rT}\n{\rS_{t}}_{\rL(\rH)}$ and $\widehat{\rC}_{\rT}=\max \{16\rM^{2}\rT, \rC_{\rT}\}$, where $\rC_{\rT}=\varepsilon^{2}\rQ^{2}\n{\rS}^{2}_{\rL^{\infty}(\rI)}\rc_{\rT}{\rm Trace }~\cQ_{\mu}$, with $\rc_{\rT}$ the positive constant of Proposition \ref{prop:key2}. Here we are using the notation  $\theta_{\rF,\rB}(s)=\theta_{\rF}(s)+\theta (s)$ with $\theta_{\rF}(s)=(\rL_{g}+\mu)s+\rQ \rS_{\infty} \theta (s)$ where $\rS_{\infty}=\n{\rS}_{\rL^{\infty}}$ and $\rL_{g}$ is a positive constant assuming that $g$ is a Lispchitz continuous function $($see Remark \ref{rem:boundedness} below$)$ and $\theta$ is the function given in \eqref{eq:theta}.
\label{prop:Gcontinuous}
\end{prop}
First, we obtain a technical estimate: 
\begin{lemma} With the notations \eqref{eq:data} or \eqref{eq:databounded} one verifies the inequality 
\begin{equation}
\EE \left [\sup_{0\le t\le \rT}\nll{\int_{0}^{t}\rS_{t-s}\big (\rB\big (\rX_{s}\big) -\rB (\big (\rY_{s}\big )\big )d\rW_{s}}^{2}_{\rH}\right ]\le \rC_{\rT}\int^{\rT}_{0}\theta \big (\nl{\rX-\rY}_{\BB_{s}}^{2}\big )ds,
\label{eq:estimatekey2refined}
\end{equation}
where $\rC_{\rT}$ is the positive constant given in Proposition \ref{prop:Gcontinuous}. 
\label{lemma:estimatekey2refined}
\end{lemma}
\par
\noindent {\sc Proof} From 	\eqref{eq:estimatekey2} we deduce 
\begin{equation}
\EE \left [\sup_{0\le t\le \rT}\nll{\int_{0}^{t}\rS_{t-s}\big (\rB\big (\rX_{s}\big) -\rB (\big (\rY_{s}\big )\big )d\rW_{s}}^{2}_{\rH}\right ]\le \rC_{\rT}\int^{\rT}_{0}\EE \big [\nl{\rB\big (\rX_{s}\big )-\rB\big (\rY_{s}\big )}^{2}_{\cL_{\cQ_{\mu}}}\big ]ds.
\label{eq:estimatekey21}
\end{equation}
On the other hand, by definition of the functional space $\cL^{2}_{\cQh}$ and \eqref{eq:moduloscontinuitybeta}, one has
$$
\begin{array}{ll}
\nl{\rB\big (\rX_{s}\big )-\rB\big (\rY_{s}\big )}^{2}_{\cL_{\cQh}}& \hspace*{-.2cm}\le \varepsilon ^{2}\rQ^{2}\n{\rS}^{2}_{\rL^{\infty}(\rI)}\big (\beta \big (\rX_{s}\big )-\beta \big (\rY_{s}\big )\big )^{2}{\rm Trace }~\cQ_{\mu} \\ [.25cm]
&  \hspace*{-.2cm}\le \varepsilon ^{2}\rQ^{2}\n{\rS}^{2}_{\rL^{\infty}(\rI)}\big (\theta \big (\rX_{s}-\rY_{s}\big )\big )^{2}{\rm Trace }~\cQ_{\mu},
\end{array}
$$
(see \eqref{eq:boundbB}), and then
$$
\begin{array}{ll}
\EE \big [\nl{\rB\big (\rX_{s}\big )-\rB\big (\rY_{s}\big )}^{2}_{\cL_{\cQ}}\big ] & 
\hspace*{-.2cm}\le \varepsilon ^{2}\rQ^{2}\n{\rS}^{2}_{\rL^{\infty}(\rI)}\EE\big [\nl{\theta \big (\rX_{s}-\rY_{s}\big )}^{2}_{\rH}\big ]{\rm Trace }~\cQ_{\mu} \\ [.25cm]
&\le \varepsilon ^{2}\rQ^{2}\n{\rS}^{2}_{\rL^{\infty}(\rI)}\theta \big (\EE\big [\nl{\rX_{s}-\rY_{s}\big )}^{2}_{\rH}\big ]\big ){\rm Trace }~\cQ_{\mu},
\end{array}
$$
(see in Section \ref{sec:coalbedofunction} the concavity and other properties of the function $\theta $). Therefore, from 
\eqref{eq:estimatekey21} the proof ends.$\fin$

\par

\smallskip

\noindent {\sc Proof of Proposition \ref{prop:Gcontinuous}.} Estimate \eqref{eq:Gcontinuous} follows from
$$
\begin{array}{ll}
\disp \EE \left [\sup_{0\le s\le t}\nl{\big (\cG^{\ru_{0}}\rX\big )_{s}-\big (\cG^{\widehat{\ru}_{0}}\rY\big )_{s}}^{2}_{\rH}\right ]  
&\disp \hspace*{-.2cm}  \le  4\sup_{0\le s\le t}\n{\rS_{s}\big (\ru_{0}-\widehat{\ru}_{0}\big )}_{\rH}^{2}\\ [.3cm]
&\disp +16t\left (\EE \left [\sup_{0\le s\le t}\int^{s}_{0}\nl{\rS_{s-\tau}\big (\rF(\rX_{\tau})-\rF(\rY_{\tau})\big )}_{\rH}^{2}d\tau\right ]\right . \\ [.3cm]
&\disp \left .+\EE\left [\sup_{0\le s\le t}\int^{s}_{0}\nl{\rS_{s-\tau}\big (\rB(\rX_{\tau})-\rB(\rY_{\tau})\big )}_{\rH}^{2}d\rW_{\tau}\right ]\right ) \\ [.5cm] 
&\disp \hspace*{-.2cm}  \le  4\rM^{2}\n{\ru_{0}-\widehat{\ru}_{0}}_{\rH}^{2} +16\rM^{2}t\EE \left [ \int^{\rT}_{0}\nl{\rF(\rX)-\rF(\rY)}_{\BB_{s}}^{2} ds \right ] \\ [.45cm] 
& \disp +\rC_{\rT}\int^{\rT}_{0}\theta \big (\nl{\rX-\rY}_{\BB_{s}}^{2}\big )ds\\ [.45cm]
&\disp \hspace*{-.2cm}  \le  4\rM^{2}\n{\ru_{0}-\widehat{\ru}_{0}}_{\rH}^{2}+16\rM^{2}t  \int^{\rT}_{0}\theta_{\rF}\big (\nl{\rX-\rY}_{\BB_{s}}^{2}\big ) ds \\ [.45cm]
& \disp +\rC_{\rT}\int^{\rT}_{0}\theta \big (\nl{\rX-\rY}_{\BB_{s}}^{2}ds  \\ [.4cm]
\end{array}
$$
(see \eqref{eq:moduloscontinuitybeta} and \eqref{eq:estimatekey2refined}).$\fin$
\par
\medskip

Now we are in a position to give the proof of the main theorem of this paper.
\par
\noindent
{\sc Proof of Theorem \ref{theo:mainTheoremintro}} As said before (see Remark \ref{rem:boundedness}),  without loss of generality we may assume that $g(u)$ is globally Lipschitz and increasing function. We may adapt the reasoning of the proof of  Proposition \ref{prop:Gcontinuous} to the functions
\begin{equation}
\left \{
\begin{array}{l}
\rF(\rX_{t})=\mu\big (\rX_{t}+10\big )-\widehat{g}\big (\rX_{t}+10\big )+\rQ \rS \theta(\rX_{t}),\\ [.1cm]
\rB(\rX_{t})=\varepsilon\rQ\rS \theta (\rX_{t}).
\end{array}
\right .
\label{eq:databounded}
\end{equation}
(see \eqref{eq:truncationg}). More precisely, since \eqref{eq:Gcontinuous} implies 
$$
\left \{
\begin{array}{ll}
\disp \nl{\cG^{\ru_{0}}\rX}_{\BB_{t}}^{2}&\hspace*{-.2cm} \le 2^{2}\left (
\nl{\cG^{\widehat{\ru}_{0}}0}_{\BB_{t}}^{2}+
\nl{\cG^{\ru_{0}}\rX-\cG^{\widehat{\ru}_{0}}0}_{\BB_{t}}^{2}
\right )\\ [.2cm]
&\hspace*{-.2cm}		
		\le 2^{2}\left (\nl{\cG^{\widehat{\ru}_{0}}0}_{\BB_{t}}^{2}+
	4\rM^{2}\n{\ru_{0}}_{\rH}^{2}+\widehat{\rC}_{\rT} \int^{\rT}_{0}\theta_{\rF,\rB}\big (\nl{\rX}_{\BB_{s}}^{2}\big )ds \right ),\quad t\in [0,\rT],
\end{array}
\right .
$$
we deduce
\begin{equation}
\n{\rX_{n+1}}_{\BB _{t}}^{2}\le \rv_{0}+\alpha \int^{t}_{0}\theta_{\rF,\rB}\big ( \n{\rX_{n}}_{\BB_{s}}^{2}\big )ds,\quad n\ge 0,
\label{eq:norsucessiveapproximatio}
\end{equation}
where $\rv_{0}= 4\big (\max\{4\rM^{2},1\}\n{\rX_{0}}_{\rH}^{2}+\nl{\cG^{\widehat{\ru}_{0}}0}_{\BB_{t}}^{2}\big )
$ and $\alpha =4\widehat{\rC}_{\rT}>0$. Next, from  Theorem \ref{theo:coalbedohybridfunction}, we consider a global solution $\rv$ on $[0,\rT]$ of the simple integral equation 
\begin{equation}
\rv(t)=\rv_{0}+\alpha  \int^{t}_{0}\theta_{\rF,\rB}\big (\rv(s)\big )ds,\quad t\in [0,\rT]. 
\label{eq:functionu1}
\end{equation}
From \eqref{eq:functionu1} we have
$$
\rv(t)-\n{\rX_{n+1}}^{2}_{\BB_{t}}\ge \alpha \int^{t}_{0}\left (\theta_{\rF,\rB}\big (\rv(s)\big )-\theta_{\rF,\rB}\big (\n{\rX_{n}}^{2}_{\BB_{s}}\big )\right )ds.
$$
Certainly $\n{\rX_{0}}^{2}_{\rH}\le \rv_{0}$ (see \eqref{eq:Picardapproximations}).  Then, the monotonicity of the function $\theta_{\rF,\rB}(\rU)$ implies, by induction, the inequality
$$
\n{\rX_{n}}^{2}_{\BB_{t}}\le \rv(t)\quad \hbox{for $t\in[0,\rT]$},
$$
where the function $\rv(t)$ is independent on $n$ (see \eqref{eq:functionu1}). 
Thus, $\{\rX_{n}\}_{n\ge 0}$ is a bounded sequence in $\BB_{\rT}$. Therefore, for each $n\ge 0$
$$
r_{n}(t)\doteq \sup_{m\ge n}\n{\rX_{m}-\rX_{n}}^{2}_{\BB_{t}}
$$
is a nonnegative, uniformly bounded, and nondecreasing function on $t\in [0,\rT]$. By construction, for each $t\in [0,\rT]$, we may consider the  nonincreasing sequence $\{r_{n}(t)\}_{n\ge 0}$. It implies the existence of a nonnegative, and nondecreasing function given by
$$
r(t)=\lim_{n\rightarrow \infty}r_{n}(t),\quad t\in [0,\rT].
$$
On the other hand, a similar reasoning as in  the proof of  Proposition \ref{prop:Gcontinuous} leads to
$$
\n{\rX_{m}-\rX_{n}}_{\BB_{t}}^{2}\le \alpha\int^{t}_{0}\theta_{\rF,\rB}\big ( \n{\rX_{m-1}(s)-\rX_{n-1}(s)}^{2}_{\BB_{s}}\big )ds.
$$
Therefore we obtain
$$
r(t)\le r_{n}(t)\le \alpha\int^{r}_{0}\theta_{\rF,\rB}\big (r_{n-1}(s)\big )ds,\quad t\in [0,\rT]
$$
whence, by the Lebesgue Convergence Theorem,
$$
r(t)\le \alpha\int^{r}_{0}\theta_{\rF,\rB}\big (r(s)\big )ds,\quad t\in [0,\rT].
$$
Now, we deduce that $r(t)\equiv 0$ for $t\in [0,\rT]$ thanks to Theorem \ref{theo:coalbedohybridfunction}. Since
$$
\n{\rX_{m}-\rX_{n}}^{p}_{\BB_{\rT}}\le r_{n}(\rT)
$$
we conclude that
$$
\n{\rX_{m}-\rX_{n}}^{2}_{\BB_{\rT}}\rightarrow 0\quad \hbox{as $m,n\rightarrow \infty$}.
$$
Therefore, the Picard type approximation, $\{\rX_{n}\}_{n\ge 0}\subset \BB_{\rT}$ given by \eqref{eq:Picardapproximations}, is a Cauchy sequence in $\BB_{\rT}$ and then we have proved the existence of a point fixed, $\rX$, of the operator $\cG^{\ru_{0}}$. 
So, we get the existence and uniqueness of the solution $\rX^{\ru_{0}}$. Moreover, the inequality \eqref{eq:Gcontinuous} becomes
\begin{equation}
\disp \n{\ru^{\ru_{0}}-\ru^{\widehat{\ru}_{0}}}_{\BB_{t}}^{2}\le
4\rM^{2}\n{\ru_{0}-\widehat{\ru}_{0}}_{\rH}^{p}+\widehat\rC_{\rT} \int^{t}_{0}\theta_{\rF,\rB} \big (\n{\ru^{\ru_{0}}-\ru^{\widehat{\ru}_{0}}}_{\BB_{s}}^{2}\big )ds,
\label{eq:Gpoinfixedeq}
\end{equation}
for $\in [0,\rT]$. Therefore, we may rewrite \eqref{eq:Gpoinfixedeq} using the integral inequality 
$$
\int^{\n{\ru^{\ru_{0}}-\ru^{\widehat{\ru}_{0}}}_{\BB_{t}}^{2}}_{4\rM^{2}\n{\ru_{0}-\widehat{\ru}_{0}}_{\rH}^{2}}\dfrac{ds}{\theta_{\rF,\rB}(s)}\le \widehat{\rC}_{\rT}t,\quad 0\le t\le \rT.
$$
\par
\noindent Finally, since we are working  in $\rH=\rL^{2}(\rI)$, we have
$$
\big (\big (\cG^{\ru_{0}}\rX\big )_{s}-\big (\cG^{\widehat{\ru}_{0}}\rX\big )_{s}\big )_{+}\le\rS_{s}\big (\ru_{0}-\widehat{\ru}_{0}\big )_{+}.
$$
Then, by the reasoning of the proof of Proposition \ref{prop:Gcontinuous}, we deduce
\begin{equation}
\disp \n{\big (\ru^{\ru_{0}}-\ru^{\widehat{\ru}_{0}}\big )_{+}}_{\BB_{t}}^{2}\le
4\rM^{2}\left (\n{\big (\ru_{0}-\widehat{\ru}_{0}\big )_{+}}_{\rH}^{2}+\widehat{\rC}_{\rT} \int^{t}_{0}\theta _{\rF,\rB}\big (\n{\big (\ru^{\ru_{0}}-\ru^{\widehat{\ru}_{0}}\big )_{+}}_{\BB_{s}}^{2}\big )ds\right ) 
\label{eq:Gpoinfixedeqcomp}
\end{equation}
for which the conclusion holds.
The proof of the comparison part of the Theorem is as follows: From \eqref{eq:Gpoinfixedeqcomp} we get
$$
\n{\big (\rX^{\ru_{0}}-\rX^{\widehat{\ru}_{0}}\big )_{+}}_{\BB_{t}}^{2}\le \widehat{\rC}_{\rT} \int^{t}_{0}\theta _{\rF,\rB}\big (\n{\big (\big (\ru^{\ru_{0}}-\ru^{\widehat{\ru}_{0}}\big )_{+}}_{\BB_{s}}^{2}\big )ds ,\quad t\in [0,\rT].
$$
Then, from the above conclusions of the Theorem we find that \ $\n{\big (\ru^{\ru_{0}}_{t}-\ru^{\widehat{\ru}_{0}}_{t}\big )_{+}}_{\rH}^{2}=0.\fin$

\par
\medskip

\begin{rem}\rm Many extensions of the main result of this paper seem to be possible. For instance, we can consider the case of function $c(x)$ taking into account the different heat capacity of continents and seas, we can replace the emitted radiative energy by a more general term of the form $g(t,x,u)$, we can consider a more general diffusion operator, and , finally, we can also consider the noise corresponding to a time-periodic Solar function $Q(t)$ (whose case deterministic was treated in~\cite{Badii}).$\fin$
\label{rem:Qperiodic}
\end{rem}
\begin{rem}\rm Other types of generalizations will be presented in the paper \cite{diaz25abstract} where, for instance, we prove the uniqueness of solutions of some abstract stochastic differential equations, under Nagumo's type conditions, which are especially useful when there are some singularities in terms that depend on time.$\fin$ 
\end{rem}

\bibliographystyle{amsplain}
\bibliography{SDEBMnonLipschitz}

\bigskip

\begin{tabular}{ll}
Gregorio Díaz & Jesús Ildefonso  Díaz \\
& Instituto Matemático Interdisciplinar (IMI)\\
Dpto. Análisis Matemático & Dpto. Análisis Matemático \\
y Matem\'atica Aplicada & y Matem\'atica Aplicada \\
U. Complutense de Madrid & U. Complutense de Madrid \\
28040 Madrid, Spain & 28040 Madrid, Spain\\
{\tt gregoriodiazdiaz@gmail.com} & {\tt jidiaz@ucm.es}
\end{tabular}

\end{document}